\documentclass[11pt]{amsart}
\usepackage{amsmath,amsthm, amscd, amssymb, amsfonts, mathtools}
\usepackage{tikz-cd}
\usepackage[all]{xy}
\usepackage{mathrsfs}
\usepackage[inline]{enumitem}
\usepackage[T2A,T1]{fontenc}
\usepackage{bigstrut}
\usepackage{rotating}
\usepackage{hyperref}
\usepackage{multirow}
\usepackage{xcolor}
\usepackage{amsmath,amsthm, amscd, amssymb, amsfonts, mathrsfs,pb-diagram,pb-xy}
\usepackage[section]{placeins}
\usepackage{multicol}
\usepackage[T1]{fontenc}
\usepackage{bm}
\usepackage{tikz}
\usepackage{float}
\usepackage{needspace}

\numberwithin{equation}{section}

\newtheorem{Theorem}{Theorem}[section]
\newtheorem*{Theorem*}{Theorem}

\newtheorem{Lemma}[Theorem]{Lemma}
\newtheorem{Proposition}[Theorem]{Proposition}
\newtheorem{Definition}[Theorem]{Definition}
\newtheorem{Example}[Theorem]{Example}

{ \theoremstyle{definition}
	 }

\newtheorem*{theorem-main}{Theorem}

\newcommand{\id}{\operatorname{id}}

\newcommand{\toba}{\mathscr B}

\newcommand{\cV}{\mathcal{V}}

\newcommand{\sele}{\mathfrak{S}}

\newcommand{\I}{{\mathbb I}}
\newcommand{\N}{{\mathbb N}}
\newcommand{\Z}{{\mathbb Z}}

\newcommand{\ot}{\otimes}
\newcommand\ad{\operatorname{ad}}

\newcommand{\eny}{\mathfrak{E}}

\newcommand{\yd}[1]{{}_{#1}^{#1}\mathcal{YD}}

\newcommand{\bq}{\mathfrak{q}}

\title[Pointed Hopf Algebras from Pale Braidings]{Pointed Hopf Algebras from Pale Braidings in Characteristic 2}

\author[Bagio,Orozco Garc\'ia]
{Dirceu Bagio, Laura Natalia Orozco Garc\'ia}

\address{D.B., L. O. G.: Departamento de Matem\'atica, Universidade Federal de Santa Catarina,
	88040-900, Florian\'opolis, SC, Brazil}
\email{d.bagio@ufsc.br,lauraorozco4e@gmail.com  }

\thanks{D. Bagio was partially supported by Funda\c c\~ao de Amparo \`a Pesquisa e Inova\c c\~ao do Estado de Santa Catarina
(FAPESC), Edital 21/2024. L. N. O. Garc\'ia was supported by CAPES-Brazil.}

\begin{document}

\begin{abstract}
	We determine Nichols algebras associated with several families of pale braidings over an algebraically closed field of characteristic $2$. For the rank-$4$ families under consideration, we obtain presentations by generators and relations, explicit monomial bases and dimension formulas. We also study families of arbitrary rank consisting of a two-dimensional pale block and several points, including cases with both vanishing and nonvanishing interaction parameters. Realizations as Yetter--Drinfeld modules over suitable finite abelian groups yield finite-dimensional pointed Hopf algebras by bosonization.
\end{abstract}

\maketitle

\section*{Introduction}

Nichols algebras are central to the study of pointed Hopf algebras. In the lifting method of Andruskiewitsch and Schneider, their determination in Yetter--Drinfeld categories is a fundamental step in the study of finite-dimensional pointed Hopf algebras; see \cite{AS}. If $V$ is a Yetter--Drinfeld module over the group algebra of a finite group $\Gamma$ and $\toba(V)$ is finite-dimensional, then the bosonization $\toba(V)\#\Bbbk\Gamma$ is a finite-dimensional pointed Hopf algebra. Explicit presentations and dimension formulas for Nichols algebras are therefore of independent interest.

Non-diagonal braidings may arise from interactions involving indecomposable Yetter--Drinfeld modules. In particular, a \emph{pale block} is indecomposable as a Yetter--Drinfeld module but decomposable as a braided vector space; see \cite{AAG}. Such modules occur even over abelian groups in positive characteristic, where the corresponding representation categories need not be semisimple. Their interactions with other homogeneous summands provide a natural setting for studying Nichols algebras beyond diagonal type.

In characteristic zero, Andruskiewitsch, Angiono and Moya Giusti classified connected rank-$4$ pale braidings whose Nichols algebras have finite Gelfand--Kirillov dimension \cite{AAG}. The analogous problem in characteristic $2$ requires separate calculations: the identity $1=-1$ affects the braidings and their defining relations, and the characteristic-zero results do not determine the corresponding Nichols algebras in characteristic $2$. Earlier computations in this characteristic include the $16$-dimensional Nichols algebra of the two-dimensional Jordan block \cite{CLW} and other finite-dimensional non-diagonal examples and their bosonizations \cite{ABFF}.

In this paper, we determine the Nichols algebras of several characteristic-$2$ families motivated by the rank-$4$ configurations in \cite{AAG}. We treat a three-dimensional pale block and a point, a two-dimensional pale block and two points; and two two-dimensional summands, including a pale block interacting with a Jordan block. We obtain defining relations, explicit monomial bases and dimensions. We also consider arbitrary-rank families consisting of a two-dimensional pale block and several points, allowing both vanishing and nonvanishing interaction parameters. The proofs use braided commutators, skew derivations and the splitting technique. Realizations over suitable finite abelian groups yield explicit finite-dimensional pointed Hopf algebras by bosonization. Our objective is to determine these families rather than to classify all Nichols algebras arising from pale braidings in characteristic $2$. \vspace{.2cm}

We summarize the results as follows.

\begin{theorem-main}
	For each braided vector space $V$ listed in Table~\ref{tab:main}, the Nichols algebra $\toba(V)$ is finite-dimensional, with dimension as indicated in the table. In each case, the cited theorem gives a presentation by generators and relations and an explicit monomial basis.
\end{theorem-main}
\vspace*{-.5cm}
\begin{table}[H]
	\centering
	\renewcommand{\arraystretch}{1.08}
	\setlength{\tabcolsep}{12pt}
	\caption{Nichols algebras from pale braidings}
	\label{tab:main}
	\begin{tabular}{lccc}
		\hline
		$V$ & $\dim V$ & $\dim\toba(V)$ & Reference \\
		\hline
		$\eny_{3,1}$ & $4$ & $2^9$ & Theorem~\ref{teo:nichols1} \\
		$\eny_{2,2}(a)$, $a\neq 0$ & $4$ & $2^6$ & Theorem~\ref{block2_plus_1point} \\
		$\eny_{2,2}$ & $4$ & $2^5$ & Theorem~\ref{block-1point-star} \\
		$\sele_{q}(1,0)$ & $4$ & $2^6$ & Theorem~\ref{two-pale-block} \\
		$\sele_{q}(1,1)$ & $4$ & $2^{12}$ & Theorem~\ref{pale-Jordan} \\
		$\eny^\star_{2,\theta-1}$ & $\theta+1$ & $2^{2\theta}$ & Theorem~\ref{block2_plus_theta-point} \\
		$\eny_{2,\theta-1}$ & $\theta+1$ & $2^{\theta+2}$ & Theorem~\ref{several-points-case-zero} \\
		$\eny_{2,\theta-1}^{\dagger}(I,J)$ & $\theta+1$ & $2^{\theta+2+|J|}$ & Theorem~\ref{several-casezeroandnot} \\
		\hline
	\end{tabular}
\end{table}

The paper is organized as follows. Section~1 recalls the necessary preliminaries on Nichols algebras, pale blocks and the splitting technique. Sections~2--4 establish the rank-$4$ results: Section~2 treats a three-dimensional pale block and a point; Section~3, a two-dimensional pale block and two points; and Section~4, two two-dimensional summands, including the pale--Jordan case. Section~5 treats the arbitrary-rank families, including the mixed case indexed by $(I,J)$. Each of Sections~2--5 also gives the corresponding realizations over finite abelian groups and pointed bosonizations.

\subsection*{Notation and conventions}
Throughout this work, $\Bbbk$ denotes an algebraically closed field of characteristic $2$, $\mathbb{N}=\{1,2,\dots\}$, and $\mathbb{N}_0=\mathbb{N}\cup\{0\}$. For $\ell,\theta\in\mathbb{N}_0$ with $\ell\leq\theta$, we write $\mathbb{I}_{\ell,\theta}:=\{\ell,\ell+1,\dots,\theta\}$ and $\mathbb{I}_\theta:=\mathbb{I}_{1,\theta}$.
For a group $\Gamma$, $\yd{\Bbbk\Gamma}$ denotes the category of Yetter--Drinfeld modules over the group algebra $\Bbbk\Gamma$. Thus, its objects are $\Gamma$-graded vector spaces $V=\bigoplus_{g\in\Gamma}V_g$ endowed with a $\Gamma$-action such that $h\cdot V_g\subseteq V_g$ for all $g,h\in\Gamma$. The Nichols algebra of $V\in\yd{\Bbbk\Gamma}$ is denoted by $\toba(V)$. Given a vector space $V$ and a subset $X\subseteq V$, we denote by $\langle X\rangle$ the vector subspace of $V$ spanned by $X$. For a Hopf algebra $H$, we denote by $G(H)=\{g\in H\,:\,g\neq 0,\,\,\Delta(g)=g\otimes g\}$ the group of group-like elements of $H$. For $g\in G(H)$, we set $P_{1,g}(H)=\{x\in H:\Delta(x)=g\otimes x+x\otimes 1\}$, the space of $(1,g)$-skew primitive elements of $H$.
	
\section{Preliminaries}

We briefly recall the basic notions that will be used throughout the paper.

\subsection{Braided vector spaces and pale blocks}\label{subsection:bvs} 
Let $\Gamma$ be an abelian group and let $\Bbbk\Gamma$ be its group algebra. The braided tensor category $\yd{\Bbbk\Gamma}$ of Yetter--Drinfeld modules over $\Bbbk\Gamma$ consists of $\Gamma$-graded vector spaces $V=\bigoplus_{g\in\Gamma} V_g$ endowed with a \(\Gamma\)-action such that $h\cdot V_g = V_g$, for all $g,h\in\Gamma$. For $V,W\in \yd{\Bbbk\Gamma}$, the braiding
$c_{V,W}:V\otimes W \to W\otimes V$ is given by
\[c_{V,W}(v\otimes w)=g\cdot w\otimes v,\qquad v\in V_g,\quad g\in\Gamma,\quad w\in W.\]

We recall that a \emph{braided vector space} is a pair $(V, c)$, where $V$ is a vector space and $c \in \operatorname{GL}(V^{\ot 2})$ is a solution of the braid equation
\begin{align*}
(c \ot \id) (\id \ot c) (c \ot \id) &= (\id \ot c) (c \ot \id)(\id \ot c).
\end{align*}
For each object $V\in \yd{\Bbbk\Gamma}$, the pair $(V,c_{V,V})$ is a braided vector space. An important class of braided vector spaces is that of diagonal type. A braided vector space $(V, c)$ is said to be of \emph{diagonal type} if there exist a basis
$(x_i)_{i\in \I_\theta }$ of $V$ and a matrix $\bq = (q_{ij})_{i,j\in \I_\theta}$
such that $q_{ij} \in \Bbbk^\times$ and $c(x_i\otimes x_j)=q_{ij}x_j\otimes x_i$ for all $i,j\in \I_\theta $. 

It is well known that if $U\in \yd{\Bbbk\Gamma}$ is indecomposable as a braided vector space then $U$ is indecomposable as a Yetter--Drinfeld module. However, the converse does not hold. In fact, there exist indecomposable $2$-dimensional Yetter--Drinfeld modules over \(\Bbbk\Gamma\) that are of diagonal type as braided vector spaces. \smallbreak

We recall Definition 3.4 from \cite{AAG}.

\begin{Definition}\label{def-paleblock} {\rm
Let $H$ be a Hopf algebra and $U\in \yd{H}$ be finite-dimensional. We say that $U$ is a pale block if it is decomposable as a braided vector space but indecomposable in $\yd{H}$.}
\end{Definition}

\begin{Example}\label{ex-pale}{\rm
Let $V=V_1\oplus V_2$ be the vector space such that $V_1=\langle x_1,x_2\rangle$, $V_2=\langle x_3\rangle$ and  let \((q_{ij})_{1\leq i,j\leq2}\) be a matrix of nonzero scalars. Consider the braiding $c$ on $V$ given by
\begin{equation}\label{eq:pale-block-point}
\bigl(c(x_i\otimes x_j)\bigr)_{1\le i,j\le 3}=\begin{pmatrix}
q_{11}x_1\otimes x_1 &q_{11}x_2\otimes x_1 &q_{12}x_3\otimes x_1\\[2mm]
q_{11}x_1\otimes x_2 &q_{11}x_2\otimes x_2 &q_{12}x_3\otimes x_2\\[2mm]
q_{21}x_1\otimes x_3 &q_{21}(x_2+x_1)\otimes x_3 &q_{22}x_3\otimes x_3\end{pmatrix}.
\end{equation}
We realize the braided vector space $V$ as an object in $\yd{\Bbbk\Gamma}$, where
$\Gamma=\mathbb Z^2$ has basis $g_1,g_2$, as follows:
\begin{align*}
 &g_1\cdot x_1=q_{11}x_1,&  &g_1\cdot x_2=q_{11}x_2,&  &g_2\cdot x_1=q_{21}x_1,& \\
 &g_2\cdot x_2=q_{21}(x_2+x_1),& &g_i\cdot x_3=q_{i2}x_3,& &i=1,2.& 
\end{align*}
The coaction is determined by $\delta(x_i)=g_1\otimes x_i$ for $i=1,2$ and $\delta(x_3)=g_2\otimes x_3$. Since \(g_2\) acts on \(V_1\) as a single Jordan block, \(V_1\) is indecomposable in \(\yd{\Bbbk\Gamma}\). On the other hand, \(V_1\) is of diagonal type as a braided vector space, and hence it is decomposable. Thus, \(V_1\) is a pale block. Hence, we say that $V$ is the sum of a pale block and a point.}
\end{Example}

\subsection{Nichols algebra}
We refer to \cite{A} for background on braided Hopf algebras and Hopf algebras in a braided tensor category. Let $H$ be a Hopf algebra and $R$ a Hopf algebra in $\yd{H}$. The braided commutator of $a,b\in R$ is $[a,b]_c=ab-(\mu\circ c)(a\otimes b)$, where $\mu$ denotes the multiplication of $R$. The braided adjoint action of $R$ is denoted by  \(\ad_c\); see \cite[p.~165]{A}. If \(a\in R\) is primitive, then $\ad_c a(b)=[a,b]_c$ for all $b\in R$.

Given \(V \in \yd{H}\), the tensor algebra \(T(V)\) is naturally a graded
Hopf algebra in \(\yd{H}\). The Nichols algebra \(\toba(V)\) is the
quotient of \(T(V)\) by the largest homogeneous Hopf ideal contained in
\(\bigoplus_{n\ge 2} T^n(V)\).

Let $\Gamma$ be an abelian group and $V\in \yd{\Bbbk\Gamma}$ with basis $(x_i)_{i\in\mathbb{I}_\theta}$. Assume that each $x_i$ is homogeneous of degree $g_i$, i.e., $\delta(x_i)=g_i\otimes x_i$. For each \(i\in\mathbb I_\theta\), there exists a skew derivation \(\partial_i\) of \(T(V)\) determined by
\[
\partial_i(x_j) = \delta_{i,j}
\quad\text{and}\quad
\partial_i(x y) = \partial_i(x)(g_i\cdot y) + x \, \partial_i(y),\]
for all \(x,y\in T(V)\) and \(i,j\in\mathbb I_\theta\). These skew derivations extend to $\toba(V)$. The following well-known fact will be used in this paper.
	
\begin{Lemma}\label{lem:skew-derivations}
Let $V\in\yd{\Bbbk\Gamma}$, and let $(x_i)_{i\in\mathbb{I}_\theta}$ be a basis of homogeneous elements of $V$. Let $(\partial_i)_{i\in\mathbb{I}_\theta}$ denote the corresponding skew
derivations. If $y\in\toba(V)$ is homogeneous of positive degree and
\[\partial_i(y)=0 \qquad \text{for all } i\in\mathbb{I}_\theta,\]
then $y=0$. 
\end{Lemma}

\subsection{The splitting technique}
Let \(V=U\oplus W\) be a direct sum of Yetter--Drinfeld modules over a
Hopf algebra \(H\). The Nichols algebra \(\toba(V)\) splits as $\toba(V)\simeq \mathcal{K}\#\toba(W)$
where $\mathcal{K}=\toba(V)^{\operatorname{co}\toba(W)}$. Moreover, by \cite[Proposition~8.6]{HS}, \(\mathcal K\) is isomorphic to the Nichols algebra of
\[
\mathcal{K}^1=\operatorname{ad}_c(\toba(W))(U).
\]
This decomposition reduces the computation of \(\toba(V)\) to that of \(\toba(\mathcal K^1)\).

\subsection{Jordan plane in characteristic $2$}
The Jordan plane in characteristic $2$ was computed in  \cite{CLW}. Let  $\epsilon\in \Bbbk^{\times}$ and $\ell \in \N_{\ge 2}$. We recall that a \emph{block}  $\cV(\epsilon,\ell)$ is a braided vector space
with a basis $(x_i)_{i\in\I_\ell}$ such that 
\begin{align}\label{equation:basis-block}
c(x_i \otimes  x_{1}) &= \epsilon x_{1} \otimes  x_i,& c(x_i \otimes  x_j) &=(\epsilon x_j+x_{j-1}) \ot  x_i,
\end{align} 
for \(i,j\in\mathbb I_\ell\), with \(j>1\).

\begin{Proposition}\label{nichols-jordanplane}
The Nichols algebra $\toba(\mathcal{V}(1,2))$ is the algebra generated by $x_1$ and $x_2$  satisfying  the following relations:
\begin{align}\label{rel-jp}
  x_1^{2} = x_2^{4}=0,&  & x_2 x_1 x_2 x_1 = x_1 x_2 x_1 x_2,&&x_2^{2} x_1 = x_1 x_2^{2} + x_1 x_2 x_1.  
\end{align}
It has dimension $16$ and basis
\begin{align}\label{bases-jp}
    \{x_1^{m_1}x_{21}^{n_1}x_2^{m_2}\,:\,m_1,n_1\in \I_{0,1},\,\,m_2\in \I_{0,3}\},
\end{align}
where $x_{21}=x_2x_1+x_1x_2$. \qed
\end{Proposition}

\subsection{Endymion algebras in characteristic $2$}\label{Endymion-algebras}
We recall a result from \cite{ABFF} concerning  the Nichols algebra associated with the pale braiding given in \eqref{eq:pale-block-point}. More precisely, when $q_{11} = q_{12}q_{21}= 1$, we denote by $\eny_{q_{12}}(q_{22})$ the braided vector space considered in Example \ref{ex-pale}.  

\begin{Proposition}[\cite{ABFF}]
 The Nichols algebra $\toba(\eny_{q_{12}}(1))$ is generated by $x_{1}$, $x_{2}$, $x_{3}$ subject to the relations:  
 \begin{align}
	x_{i}^{2}&=0, && i\in \I_{3}, \\
    x_{1}x_{2}&=x_{2}x_{1},  && x_{1}x_{3}=q_{12}x_{3}x_{1},\\
     x_{23}^{2}&=0,          &&\text{where }\, x_{23}:=x_{2}x_{3}+ q_{12}x_{3}x_{2}.
\end{align}
Moreover, the set of monomials $\{x_{1}^{m_1}x_{2}^{m_2}x_{23}^{n_1}x_{3}^{m_3}:m_1,m_2,m_3,n_1\in\I_{0,1}\}$ is a basis of $\toba(\eny_{q_{12}}(1))$. Hence $\dim\toba(\eny_{q_{12}}(1))=2^4.$  \qed
\end{Proposition}

\section{A 3-dimensional pale block and a point}\label{sec2}
In this section, we determine the Nichols algebra associated with a braided vector space consisting of a \(3\)-dimensional pale block and a point. We use this Nichols algebra to construct an example of a finite-dimensional pointed Hopf algebra over a field of characteristic \(2\).

\subsection{Nichols algebra} Let  $V = V_{g_1} \oplus V_{g_2}$ be the braided vector space such that $V_{g_1}$ is a pale block of eigenvalue $1$ with basis 
\((x_i)_{i\in\I_3}\) and $V_{g_2} = \langle x_4 \rangle$. Explicitly, the braiding on $V$ is given by
\begin{align*}
c(x_i\otimes x_j)&=x_j\otimes x_i,& c(x_i\otimes x_4)&=q_{12}\,x_4\otimes x_i,\\
c(x_4\otimes x_i)&=q_{21}(x_i+x_{i-1})\otimes x_4,&c(x_4\otimes x_4)&=x_4\otimes x_4,&&   i,j\in\I_3,
\end{align*}
where $x_0=0$ and $q_{12} q_{21} = 1$. Let $\Gamma:=\langle g_1,g_2\rangle\simeq\Z^2$. We realize $V$ as a Yetter--Drinfeld module over $\Bbbk\Gamma$ with action given by
\begin{align}\label{action-pale3-1point}
\begin{aligned}
	g_1 \cdot x_i &= x_i, \qquad\quad &   g_2 \cdot x_i &= q_{21} (x_i + x_{i-1}), \quad i \in \I_3, \\
	  g_1 \cdot x_4 &= q_{12} x_4,\qquad\quad& g_2 \cdot x_4 &= x_4.
\end{aligned}
\end{align}
The comodule structure on \(V\) is given by
\begin{align}\label{comodule-pale-1point}
  \delta(x_i)=g_1\otimes x_i,\,\,i\in \I_3,& &\delta(x_4)=g_2\otimes x_4.   
\end{align}
We denote this braided vector space by $\eny_{3,1}$. For \(i\in\{2,3\}\), set
\[x_{4i} := x_4 x_i + q_{21} (x_i + x_{i-1}) x_4, \quad   w := x_{42} x_3 + q_{21} (x_3 + x_2) x_{42}.\]
The following theorem gives a presentation of the Nichols algebra $\toba(\eny_{3,1})$.

\begin{Theorem}\label{teo:nichols1}
	$\toba(\eny_{3,1})$ is generated by $x_1, x_2, x_3, x_4$, subject to the following relations:
	\begin{align}
		x_i^2 &= 0, \quad i \in \I_3,&  x_i x_j &= x_j x_i, \quad  i \neq j \in \I_3, \label{1point-eq1} \\[.2em]
		x_4^2 &= 0,  &x_1 x_4 &= q_{12} x_4 x_1, \label{1point-eq2} \\[.2em]
        x_{24}&:= x_2x_4+q_{12}x_{4}x_{2},  &    x_{24}^2&=0, \label{1point-eq3}  \\[.2em]
    x_{43} x_{42} x_{43} x_{42} &+ x_{42} x_{43} x_{42} x_{43}=0, &  x_{43}^{\,2} x_{42}& + x_{42} x_{43}^{\,2} + x_{42}x_{43}x_{42} = 0, \label{1point-eq5}\\[.2em]
        x_{42}^2&=0, & x_{43}^{4}&=0, \label{1point-eq6} \\[.2em]
		w^2&=0,  &    x_{43} w&= q_{21} w x_{43}, \label{1point-eq4}
	\end{align}
A basis of $\toba(\eny_{3,1})$ is given by the monomials
\begin{align}\label{monomios1}
x_1^{m_1}x_2^{m_2}x_3^{m_3}w^{l_1}x_{43}^{n_2}\,z_{23}^{l_2}\, x_{42}^{\,n_1}\, x_4^{m_4}
\end{align}
where $z_{23}=x_{42}x_{43}+x_{43}x_{42}$, $m_i \in \I_{0,1}$, $n_1,\,l_j \in \I_{0,1}$ and $n_2 \in \I_{0,3}$.  
Hence, $\dim \toba(\eny_{3,1}) = 2^9$.
\end{Theorem}

\begin{proof}
The relations \eqref{1point-eq1}, \eqref{1point-eq2}, \eqref{1point-eq3} hold in $\toba(\eny_{3,1})$ because the braiding on $\langle x_{1},x_{2},x_{3}\rangle$ is  the flip and $\langle x_{1},x_{2},x_{4}\rangle\simeq \eny_{q_{12}}(1)$ as braided vector spaces, where $\eny_{q_{12}}(1)$ is the braided vector space defined in Subsection \ref{Endymion-algebras}. 
Next, let $\mathcal{K}^{1}=\ad_c(\toba(V_{g_2})) (V_{g_1})$. For \(i\in\{2,3\}\), the coaction of \(x_{4i}\) is given by $\delta(x_{4i})=x_{4}g_{1}\otimes x_{i-1}+g_{1}g_{2}\otimes x_{4i}$. Hence, the braiding on $\mathcal{K}^{1}$ is given by
\begin{align*}
 (c(x_{4i}\otimes x_{4j}))_{i,j\in \I_{2,3}}=\begin{pmatrix}
x_{42}\otimes  x_{42} & (x_{43}+x_{42})\otimes x_{42}  \\
 x_{42}\otimes x_{43}& (x_{43}+x_{42})\otimes  x_{43}  \\
\end{pmatrix}.
\end{align*}
Therefore, by Proposition \ref{nichols-jordanplane}, the relations \eqref{1point-eq5} and \eqref{1point-eq6} hold. The relations \eqref{1point-eq1}--\eqref{1point-eq6}, together with the definitions of
$x_{42}$, $x_{43}$, $w$, and $z_{23}$, imply the following identities:
\begin{equation}\label{eq:conjunto}
\begin{aligned}
    x_{1}x_{4i}&=q_{12}x_{4i}x_{1},  \text{ for }i=2,3, & \quad  x_{42}x_3&=w+q_{21}(x_3+x_2)x_{42},\\
    x_{42}x_2&=q_{21}(x_{2}+x_{1})x_{42},   & \quad  x_{4}x_{42}&= q_{21}x_{42}x_{4},\\
    x_{42}w&=q_{21}wx_{42} ,  & \quad  x_{43}x_2&=w+q_{21}(x_{2}+x_{1})x_{43}, \\
    x_{43}x_3&=q_{21}(x_{3}+x_{2})x_{43},   & \quad  x_{4}x_{43}&=q_{21}(x_{43}+x_{42})x_{4} \\
    wx_2&=q_{21}(x_{2}+x_{1})w,   & \quad  wx_3&=q_{21}(x_{3}+x_{2})w, \\
    x_{4}w&=q_{21}^{2}wx_{4}+q_{21}z_{23} & \quad x_{1}w&=q_{12}wx_{1}.
\end{aligned}
\end{equation}
To verify \eqref{1point-eq4}, we use the skew derivations, which satisfy

\[
\begin{array}{c|ccc}
 & x_{42} & x_{43} & w \\ \hline
\partial_1 & x_4 & 0   & x_{43} \\
\partial_2 & 0   & x_4 & x_{42} \\
\partial_3 & 0   & 0   & 0      \\
\partial_4 & 0   & 0   & 0
\end{array}
\]
Using these identities, we obtain
\[\partial_i(w^2)=0,\qquad
\partial_i(x_{43}w+q_{21}wx_{43})=0,
\qquad i=3,4.\]
Also, $\partial_{1}(x_{43}w+q_{21}wx_{43})=0$ and $\partial_{2}(x_{43}w+q_{21}wx_{43})=q_{12}x_4w+x_{43}x_{42}+x_{42}x_{43}+q_{21}wx_4=0$ by \eqref{eq:conjunto}. Hence, the relation $x_{43}w=q_{21}wx_{43}$ holds in $\toba(\eny_{3,1})$. Using these relations, we compute $\partial_{1}(w^{2})=q_{12}x_{43}w+wx_{43}=0$, and similarly $\partial_2(w^{2})=q_{12}x_{42}w+wx_{42}=0$ by \eqref{eq:conjunto}. Consequently, \eqref{1point-eq4} holds in $\toba(\eny_{3,1})$. 

Let $\toba$ be the algebra defined by the presentation in the statement
of the theorem. By the previous paragraph, there is a well-defined
algebra morphism
\[
\pi:\toba\longrightarrow\toba(\eny_{3,1}),
\qquad x_i\longmapsto x_i,\quad i\in\I_4.
\]
Since $\toba(\eny_{3,1})$ is generated by $x_1,x_2,x_3,x_4$, the map
$\pi$ is surjective. Using the identities in \eqref{eq:conjunto} together with the defining relations, we obtain that every monomial in
$x_1$, $x_2$, $x_3$ and $x_4$ can be written as a linear combination of the
monomials in \eqref{monomios1}. Hence, these monomials span $\toba$,
and, since $\pi$ is surjective, their images span
$\toba(\eny_{3,1})$.

We now prove that the monomials in \eqref{monomios1} are linearly
independent. Suppose, to the contrary, that there exists a nonzero
homogeneous linear combination $S$ of these monomials which vanishes,
and choose such an $S$ of minimal degree. Since
\[\partial_4
\left(
x_1^{m_1}x_2^{m_2}x_3^{m_3}w^{l_1}
x_{43}^{n_2}z_{23}^{l_2}x_{42}^{n_1}x_4^{m_4}
\right)
=
\delta_{1,m_4}
x_1^{m_1}x_2^{m_2}x_3^{m_3}w^{l_1}
x_{43}^{n_2}z_{23}^{l_2}x_{42}^{n_1},\]
the minimality of $S$ implies that $m_4=0$ for every monomial occurring
in $S$ with nonzero coefficient. Set $D=\partial_4\partial_1$ and $E=\partial_4\partial_2$. We have
\begin{align*}
D\left(
x_1^{m_1}x_2^{m_2}x_3^{m_3}w^{l_1}
x_{43}^{n_2}z_{23}^{l_2}x_{42}
\right)
&=
x_1^{m_1}x_2^{m_2}x_3^{m_3}w^{l_1}
x_{43}^{n_2}z_{23}^{l_2},\\
D^2\left(
x_1^{m_1}x_2^{m_2}x_3^{m_3}w^{l_1}
x_{43}^{n_2}z_{23}
\right)
&=
x_1^{m_1}x_2^{m_2}x_3^{m_3}w^{l_1}
x_{43}^{n_2},\\
D^{n_2-1}E\left(
x_1^{m_1}x_2^{m_2}x_3^{m_3}w^{l_1}x_{43}^{n_2}
\right)
&=
x_1^{m_1}x_2^{m_2}x_3^{m_3}w^{l_1},
\qquad n_2\geq1,\\
\partial_4\partial_1\partial_2\left(
x_1^{m_1}x_2^{m_2}x_3^{m_3}w
\right)
&=
x_1^{m_1}x_2^{m_2}x_3^{m_3},
\end{align*}
which implies that $n_1=l_2=n_2=l_1=0$.  Hence $S$ is a nonzero linear
combination of
\[
x_1^{m_1}x_2^{m_2}x_3^{m_3},
\qquad m_i\in\{0,1\},
\]
which are linearly independent, a contradiction. Therefore the
monomials in \eqref{monomios1} are linearly independent in
$\toba(\eny_{3,1})$. Since they span $\toba(\eny_{3,1})$, they form a
basis. Consequently, $\toba\simeq\toba(\eny_{3,1})$.
\end{proof}

\subsection{A realization of $\toba(\eny_{3,1})$}
Let $q_{12}$ be a primitive $n$-th root of unity and 
\[
\widetilde{\Gamma}:=C_n\times C_m,
\qquad
m:=\operatorname{lcm}(4,n)\overset{(*)}{=}4n,
\]
where $(*)$ follows since $n$ is odd, as $\operatorname{char}\Bbbk=2$.
Then $\eny_{3,1}\in\yd{\Bbbk\widetilde{\Gamma}}$, with action and coaction
as in the case of $\Bbbk\Gamma$ considered in the previous subsection.
Here, we set $x_0:=0$.

\begin{Proposition}
Let $H=\toba(\eny_{3,1})\#\Bbbk\widetilde{\Gamma}$ and $\widetilde{\Gamma}=\langle g_1\rangle\times\langle g_2\rangle$. Then $H$ is a finite-dimensional pointed Hopf algebra generated by
$x_1,x_2,x_3,x_4,g_1,g_2$ satisfying \eqref{1point-eq1}--\eqref{1point-eq4} and the following relations: 
\begin{align*}
g_1x_i&=x_ig_1,\quad i\in\I_3,&g_1x_4&=q_{12}x_4g_1,\\
g_2x_i&=q_{21}(x_i+x_{i-1})g_2, \quad i\in\I_3,&g_2x_4&=x_4g_2,\\
g_1^n&=g_2^m=1,& g_1g_2&=g_2g_1.
\end{align*}
The comultiplication is determined by $g_i\in G(H)$, $i\in\I_2$, $x_j\in P_{1,g_1}(H)$, $j\in\I_3$, and $x_4\in P_{1,g_2}(H)$.
A basis of $H$ is given by the monomials
\[x_1^{m_1}x_2^{m_2}x_3^{m_3}w^{l_1}x_{43}^{n_2}z_{23}^{l_2}x_{42}^{n_1}x_4^{m_4}g_1^{s_1}g_2^{s_2},\]
where $m_i,n_1,l_j\in\{0,1\}$, $n_2\in\I_{0,3}$, $s_1\in\I_{0,n-1}$ and $s_2\in\I_{0,m-1}$. Hence, $\dim H=2^9nm$. \qed 
\end{Proposition}

\section{A pale block of dimension $2$ and two points}\label{sec-pale-2-points}

In this section, we study the Nichols algebra associated with a Yetter--Drinfeld module consisting of a pale block of dimension~$2$ and two points. 

Consider the Yetter--Drinfeld module 
$V = V_{g_1} \oplus V_{g_2} \oplus V_{g_3}$ over $\Bbbk\Gamma$, where $\Gamma=\langle g_1\rangle\times \langle g_2\rangle\times\langle g_3\rangle\simeq\Z^3$. The component $V_{g_1}=\langle x_1,x_{\frac32}\rangle$ is a pale block with eigenvalue $1$, while $V_{g_2}=\langle x_2\rangle$ and $V_{g_3}=\langle x_3\rangle$. 
Explicitly, 
the action is given by
\begin{align}\label{action-pale-2point}
\begin{aligned}
g_1 \cdot x_1 &= x_1, &
g_2 \cdot x_1 &= q_{21}x_1, &
g_3 \cdot x_1 &= q_{31}x_1,\\[1pt]
g_1 \cdot x_{\frac{3}{2}} &= x_{\frac{3}{2}}, &
g_2 \cdot x_{\frac{3}{2}} &= q_{21}(x_{\frac{3}{2}}+x_1), &
g_3 \cdot x_{\frac{3}{2}}&= q_{31}(x_{\frac{3}{2}}+a x_1), \\[1pt]
g_1 \cdot x_2 &= q_{12}x_2, &
g_2 \cdot x_2 &=x_2, &
g_3 \cdot x_2 &= q_{32}x_2, \\[1pt]
g_1 \cdot x_3 &= q_{13}x_3, &
g_2 \cdot x_3 &= q_{23}x_3, &
g_3 \cdot x_3 &=x_3,
\end{aligned}
\end{align}
where $q_{ij}\in \Bbbk^\times$. We assume that $
q_{12}q_{21}=q_{13}q_{31}=q_{23}q_{32}=1$ and $a\in \Bbbk$. The comodule structure of $V$ is given by
\begin{align}\label{coaction-pale-2point}
\delta(x_i)&=g_i\otimes x_i, \quad i\in \I_3, & \delta(x_{\frac 32})&=g_1\otimes x_{\frac 32}.
\end{align}
The previous Yetter--Drinfeld module $V$ will be denoted by $\eny_{2,2}(a)$. Throughout this section, we set
\[
x_{\frac32 2}:=x_{\frac32}x_2+q_{12}x_2x_{\frac32},
\qquad
x_{\frac32 3}:=x_{\frac32}x_3+q_{13}x_3x_{\frac32}.
\]

\subsection{Nichols algebra of $\eny_{2,2}(a)$, $a\neq 0$}
Assume that $a\in \Bbbk^\times$. The Nichols algebra $\toba(\eny_{2,2}(a))$ is presented by
generators and relations in the next result.

\begin{Theorem} \label{block2_plus_1point}
The Nichols algebra $\toba(\eny_{2,2}(a))$ is generated by
$x_1$, $x_{\frac 32}$, $x_2$, and $x_3$, with defining relations 
\begin{align}
x_1^{2} &= 0, & x_{\frac 32}^{2} &= 0, & x_{2}^{2}&=0, \label{eq1-end-211}\\
x_1x_{\frac 32} &= x_{\frac 32}x_1,& x_1x_2 &= q_{12}x_2x_1,  & x_{\frac 32 2}^{2} &= 0,\label{eq2-end-211}  \\
x_{\frac 32 3}^{2} &= 0,  &  x_3^{2} &= 0, & x_1 x_3 &= q_{13} x_3 x_1,  \label{eq4-end-211}\\
  x_2 x_3 &= q_{23}\,x_{3}x_{2}, &x_3 x_{\frac 32 2} &= q_{32}q_{31} x_{\frac 32 2} x_3\label{eq3-end-211}.
\end{align}
A basis of $\toba(\eny_{2,2}(a))$ is given by the monomials
\begin{align}\label{monomials-end}
x_1^{m_1}\, x_{\frac{3}{2}}^{n_1}\,
x_2^{m_2}\, x_{\frac{3}{2}2}^{l_2}\,
x_3^{m_3}\, x_{\frac{3}{2}3}^{l_3},
\end{align}
where $m_i,\, n_1,\, l_2,\, l_3 \in \I_{0,1}$. Hence $\dim \toba(\eny_{2,2}(a))=2^6$.
\end{Theorem}

\begin{proof}
   The relations \eqref{eq1-end-211}, \eqref{eq2-end-211}, \eqref{eq4-end-211} hold in $\toba(\eny_{2,2}(a))$ since $\langle x_{1},x_{\frac 32},x_{2}\rangle \simeq \eny_{q_{12}}(1)$ and $\langle x_{1},x_{\frac 32},x_{3}\rangle\simeq \eny_{q_{13}}(1)$ as braided vector spaces; see Subsection \ref{Endymion-algebras}. The relations
in \eqref{eq3-end-211} follow by a direct computation with the skew derivations.

Let $\toba$ be the algebra defined by the presentation above. By the preceding paragraph, there is a surjective algebra morphism  $\pi:\toba \to \toba(\eny_{2,2}(a))$. The defining relations further imply
   \begin{equation*}
\begin{aligned}
    x_1x_{\frac 32 2} &= q_{12}\, x_{\frac 32 2}x_{1},& x_1x_{\frac 32 3} &= q_{13}\,x_{\frac 32 3}x_1,  & x_{\frac 32}x_{\frac 32 2} &= q_{12}\,x_{\frac 32 2}x_{\frac 32}, \\
    x_{\frac 32}x_{\frac 32 3}&=q_{13}x_{\frac 32 3}x_{\frac 32}, & x_{\frac 32 2} x_{\frac 32 3} &= q_{23}\,q_{13}\,q_{21}\, x_{\frac 32 3} x_{\frac 32 2}, & x_2x_{\frac 32 3} &= q_{23}\, q_{21} x_{\frac 32 3}x_2.  
\end{aligned}
\end{equation*}
Together with the definitions of $x_{\frac32 2}$ and $x_{\frac32 3}$, these relations show that the monomials in \eqref{monomials-end} span $\toba$, and hence also
$\toba(\eny_{2,2}(a))$. It remains to prove their linear independence. The monomials $x_1^{m_1}x_{\frac 32}^{n_1}x_{2}^{m_2}x_{\frac 32 2}^{l_2}$ and $x_{3}^{m_3}x_{\frac 32 3}^{l_3}$ are linearly independent in $\toba(\eny_{2,2}(a))$, as follows from the embeddings 
   $$\toba(\eny_{q_{12}}(1))\simeq \Bbbk \langle x_{1},x_{\frac 32},x_{2}\rangle \hookrightarrow \toba(\eny_{2,2}(a)) \hookleftarrow \Bbbk \langle x_{1},x_{\frac 32}, x_3 \rangle \simeq \toba(\eny_{q_{13}}(1)).$$
Finally, the decomposition $V=(V_{g_1}\oplus V_{g_2})\oplus V_{g_3}$ and the splitting technique yield a linear isomorphism $\toba(\eny_{2,2}(a))\simeq\toba(\eny_{q_{12}}(1))\otimes\mathcal K$, where $\mathcal K^1=\ad_c\bigl(\toba(\eny_{q_{12}}(1))\bigr)(V_{g_3})$ and $\mathcal K=\toba(\mathcal K^1)$. Since $x_3,x_{\frac32 3}\in\mathcal K$, the preceding linear independence and the above factorization imply that the monomials in
\eqref{monomials-end} are linearly independent. Therefore, they form a basis of $\toba(\eny_{2,2}(a))$.
\end{proof}

\subsubsection{A realization of $\toba(\eny_{2,2}(a))$}
Let $q_{12}$, $q_{13}$, and $q_{23}$ be primitive $n_1$-th, $n_2$-th and $n_3$-th roots of unity, respectively, and let
\[\widetilde{\Gamma}=C_{m_1}\times C_{m_2}\times C_{m_3},\]
where $m_1=\operatorname{lcm}(n_1,n_2)$, $m_2=\operatorname{lcm}(2,n_1,n_3)$ and $m_3=\operatorname{lcm}(2,n_2,n_3)$.
Then $\eny_{2,2}(a)\in\yd{\Bbbk\widetilde{\Gamma}}$, with action and coaction as in the case of $\Bbbk\Gamma$ considered in the previous
subsection.

\begin{Proposition}
Let $H=\toba(\eny_{2,2}(a))\#\Bbbk\widetilde{\Gamma}$ and $\widetilde{\Gamma}=\langle g_1\rangle\times\langle g_2\rangle\times\langle g_3\rangle$. Then $H$ is a finite-dimensional pointed Hopf algebra generated by $x_1$,$\,x_{\frac32}$,$\,x_2$,$\,x_3$,$\,g_1$,\newline$\,g_2$ and $g_3$ satisfying
\eqref{eq1-end-211}--\eqref{eq3-end-211} and the following relations:
\begin{align*}
g_1x_i&=q_{1i}x_ig_1,\quad i\in\I_3,&g_1x_{\frac32}&=x_{\frac32}g_1,\\
g_2x_i&=q_{2i}x_ig_2,\quad i\in\I_3,&g_2x_{\frac32}&=q_{21}(x_{\frac32}+x_1)g_2,\\
g_3x_i&=q_{3i}x_ig_3,\quad i\in\I_3,&g_3x_{\frac32}&=q_{31}(x_{\frac32}+ax_1)g_3,\\
g_1^{m_1}&=g_2^{m_2}=g_3^{m_3}=1,&g_ig_j&=g_jg_i,\quad i,j\in\I_3.\end{align*}
The comultiplication is determined by $g_i\in G(H)$, $i\in\I_3$, $x_i\in P_{1,g_i}(H)$, $i\in\I_3$ and $x_{\frac32}\in P_{1,g_1}(H)$.
A basis of $H$ is given by the monomials
\[x_1^{r_1}x_{\frac32}^{r_{\frac32}}x_2^{r_2}x_{\frac32 2}^{l_2}x_3^{r_3}x_{\frac32 3}^{l_3}g_1^{s_1}g_2^{s_2}g_3^{s_3},\]
where $r_1,r_{\frac32},r_2,r_3,l_2,l_3\in\I_{0,1}$, $s_i\in\I_{0,m_i-1}$ and $i\in\I_3$. Hence $\dim H=2^6m_1m_2m_3$. \qed 
\end{Proposition}

\subsection{Nichols algebra in the case $a=0$}
In this subsection, we write $\eny_{2,2}$ for the braided vector space
$\eny_{2,2}(0)$.

\begin{Theorem}\label{block-1point-star}
The Nichols algebra $\toba(\eny_{2,2})$ is generated by
$x_1, x_{\frac 32}, x_2, x_3$ with defining relations
\begin{align}
x_{\frac 32}x_1&= x_1x_{\frac 32}, &x_{i}^2 &=0,     & i\in&\,\,\I_2, \label{eq1-211-new}\\
x_{1}x_2 &= q_{12}x_2 x_{1}       , & x_{\frac 32 2}^2&=0,     &x_{\frac 32}^2&=0,\label{eq2-211-new}\\
x_{\frac 32} x_3 &= q_{13} x_3 x_{\frac 32 },   & x_{1}x_3&=q_{13}x_3x_1, & x_{3}^2&=0, \label{eq3-211-new}\\
x_2x_3&= q_{23}x_3x_2\label{eq4-211-new}. 
\end{align}
A basis of $\toba(\eny_{2,2})$ is given by the monomials 
\begin{align}\label{monomials-end211}
x_{2}^{m_2}x_{\frac 32 2}^{l_1}\,x_3^{m_3}\,x_{\frac 32}^{n_1}x_1^{m_1},
\end{align}
where $m_{i},n_1,l_1\in \{0,1\}$ with $i\in \I_{3}$. Hence $\dim \toba(\eny_{2,2}) = 2^{5}$.
\end{Theorem}

\begin{proof}
Observe that \eqref{eq1-211-new} and \eqref{eq2-211-new} hold in $\toba(\eny_{2,2})$, since $\langle x_1,x_{\frac 32 }, x_2 \rangle\simeq \eny_{q_{12}}(1)$ as a braided vector space. On the other hand, the relations \eqref{eq3-211-new} and \eqref{eq4-211-new} are easily verified via skew derivations.

Let $\toba$ be the algebra defined by the presentation above. By the preceding paragraph, there is a surjective algebra morphism $\pi:\toba\longrightarrow\toba(\eny_{2,2})$. The defining relations imply
\begin{equation}\label{eq:conjunto-end211}
\begin{aligned}
x_1x_{\frac32 2}    &=q_{12}x_{\frac32 2}x_1,&x_{\frac32}x_{\frac32 2}    &=q_{12}x_{\frac32 2}x_{\frac32},\\
x_2x_{\frac32 2}    &=q_{21}x_{\frac32 2}x_2,&x_3x_{\frac32 2}    &=q_{31}q_{32}x_{\frac32 2}x_3.
\end{aligned}
\end{equation}
Together with the definition of $x_{\frac32 2}$, these relations show that the monomials in \eqref{monomials-end211} span $\toba$, and hence also
$\toba(\eny_{2,2})$.

It remains to prove that the monomials in \eqref{monomials-end211} are linearly independent. Suppose, to the contrary, that there exists a
nontrivial homogeneous linear combination $S$ of these monomials that vanishes, and choose such an $S$ of minimal degree.
Since
\begin{align*}
\partial_1\bigl(x_2^{m_2}x_{\frac32 2}^{l_1}x_3^{m_3}x_{\frac32}^{n_1}x_1^{m_1}\bigr)&=\delta_{m_1,1}x_2^{m_2}x_{\frac32 2}^{l_1}x_3^{m_3}x_{\frac32}^{n_1},\\
\partial_{\frac32}\bigl(x_2^{m_2}x_{\frac32 2}^{l_1}x_3^{m_3}x_{\frac32}^{n_1}\bigr)&=\delta_{n_1,1}x_2^{m_2}x_{\frac32 2}^{l_1}x_3^{m_3},\\
\partial_3\bigl(x_2^{m_2}x_{\frac32 2}^{l_1}x_3^{m_3}\bigr)&=\delta_{m_3,1}x_2^{m_2}x_{\frac32 2}^{l_1},\\
\partial_1\partial_2\bigl(x_2^{m_2}x_{\frac32 2}^{l_1}\bigr)&=\delta_{l_1,1}x_2^{m_2},\\\partial_2(x_2^{m_2})&=\delta_{m_2,1},\end{align*}
the minimality of $S$ successively forces $m_1=n_1=m_3=l_1=m_2=0$ for every monomial occurring in $S$, a contradiction.
Therefore, these monomials are linearly independent in $\toba(\eny_{2,2})$, and consequently $\toba \simeq \toba(\eny_{2,2})$.
    \end{proof}

\subsubsection{A realization of $\toba(\eny_{2,2})$}
Let $q_{12}$, $q_{13}$, and $q_{23}$ be primitive $n_1$-th, $n_2$-th and $n_3$-th roots of unity, respectively, and let
\[
\widetilde{\Gamma}=C_{m_1}\times C_{m_2}\times C_{m_3},\]
where $m_1=\operatorname{lcm}(n_1,n_2)$, $m_2=\operatorname{lcm}(2,n_1,n_3)$ and $m_3=\operatorname{lcm}(n_2,n_3)$.
Then $\eny_{2,2}\in\yd{\Bbbk\widetilde{\Gamma}}$, with action and coaction
as in the case of $\Bbbk\Gamma$ considered in the previous subsection.

\begin{Proposition}
Let $H=\toba(\eny_{2,2})\#\Bbbk\widetilde{\Gamma}$ and
$\widetilde{\Gamma}
=\langle g_1\rangle\times\langle g_2\rangle\times\langle g_3\rangle$.
Then $H$ is a finite-dimensional pointed Hopf algebra generated by
$x_1$, $x_{\frac32}$, $x_2$, $x_3$, $g_1$, $g_2$ and $g_3$ satisfying
\eqref{eq1-211-new}--\eqref{eq4-211-new} and the following relations:
\begin{align*}
g_1x_i&=q_{1i}x_ig_1,\quad i\in\I_3,&g_1x_{\frac32}&=x_{\frac32}g_1,\\
g_2x_i&=q_{2i}x_ig_2,\quad i\in\I_3,&g_2x_{\frac32}&=q_{21}(x_{\frac32}+x_1)g_2,\\
g_3x_i&=q_{3i}x_ig_3,\quad i\in\I_3,&g_3x_{\frac32}&=q_{31}x_{\frac32}g_3,\\
g_1^{m_1}&=g_2^{m_2}=g_3^{m_3}=1,&g_ig_j&=g_jg_i,\quad i,j\in\I_3.\end{align*}
The comultiplication is determined by $g_i\in G(H)$, $i\in\I_3$,
$x_i\in P_{1,g_i}(H)$, $i\in\I_3$, and
$x_{\frac32}\in P_{1,g_1}(H)$.
A basis of $H$ is given by the monomials
\[
x_2^{r_2}x_{\frac32 2}^{l_1}x_3^{r_3}
x_{\frac32}^{r_{\frac32}}x_1^{r_1}
g_1^{s_1}g_2^{s_2}g_3^{s_3},
\]
where $r_1,r_{\frac32},r_2,r_3,l_1\in\I_{0,1}$ and $s_i\in\I_{0,m_i-1},\quad i\in\I_3$. Hence $\dim H=2^5m_1m_2m_3$.
\end{Proposition}

\section{Two blocks of dimension $2$} 
We next study braided vector spaces consisting of two blocks. In this section, $\Gamma=\Z^2$ denotes the abelian group generated by $g_1$ and $g_2$ and  $V = V_{g_1} \oplus V_{g_2}$ where 
$V_{g_1} = \langle x_1, x_{\frac{3}{2}} \rangle$ and
$V_{g_2} = \langle x_2, x_{\frac{5}{2}} \rangle$. We realize $V$ as a Yetter--Drinfeld module over $\Bbbk\Gamma$ with action given by
\begin{equation}\label{action-two-blocks}
\begin{aligned}
g_1 \cdot x_1 &= x_1, 
&\qquad g_2 \cdot x_1 &= q_{21} x_1, \\
g_1 \cdot x_{\frac{3}{2}} &= x_{\frac{3}{2}}, 
&\qquad g_2 \cdot x_{\frac{3}{2}} &= q_{21} (x_1 +  x_{\frac{3}{2}}), \\
g_1 \cdot x_2 &= q_{12} x_2, 
&\qquad g_2 \cdot x_2 &= x_2, \\
g_1 \cdot x_{\frac{5}{2}} &= q_{12} (x_{\frac{5}{2}} + a x_2), 
&\qquad g_2 \cdot x_{\frac{5}{2}} &= x_{\frac{5}{2}} + bx_2.
\end{aligned}
\end{equation}
The coaction is determined by
\begin{align*}
    &\delta(x_i)=g_i\otimes x_i,& &\delta(x_{i+\frac{1}{2}})=g_i\otimes x_{i+\frac{1}{2}},& &i\in \I_2.&
\end{align*}
The previous braiding depends on the parameters \(a,b\in\Bbbk\). In what follows in this section, we consider the cases $(a,b)=(1,0)$ and $(a,b)=(1,1)$, which correspond to the characteristic-$2$ counterparts of some of the cases studied in \cite[Section~6]{AAG}. Notice that, in characteristic $2$, the distinction between $1$ and $-1$ disappears; consequently, some of the cases considered in \cite{AAG} coincide in the present setting.

We assume that $q_{12}=q=q_{21}^{-1}$ and we denote this braided vector space by $\sele_{q}(a,b)$. Throughout this section, we use the notation
\begin{align*}
x_{\frac32 2}&:=x_{\frac32}x_2+q_{12}x_2x_{\frac32}, & x_{\frac32\frac52}&:=x_{\frac32}x_{\frac52}+q_{12}\big(x_{\frac52}+ax_2\big)x_{\frac32},\\
x_{2\frac52}&:=x_2x_{\frac52}+\big(x_{\frac52}+bx_2\big) x_2,& x_{1{\frac{5}{2}}}&:=x_{1}x_{\frac{5}{2}}+q_{12}(x_{\frac{5}{2}}+ax_2).
\end{align*}

\subsection{Nichols algebra of two pale blocks}\label{subsec41} In this subsection, we assume that $a=1$ and $b=0$. The following theorem describes this case.

\begin{Theorem}\label{two-pale-block}
The Nichols algebra $\toba(\sele_{q}(1,0))$ is generated by 
$x_1, x_{\frac{3}{2}}, x_2, x_{\frac{5}{2}}$ with defining relations
\begin{align}
 x_i^2                &=0,  \quad i\in \I_2,       & x_{i+\frac{1}{2}}^2 &= 0,\quad i\in \I_2,    \label{eq:5.1}\\
 x_i x_{i+\frac{1}{2}} &=  x_{i+\frac{1}{2}} x_i,  \quad i\in \I_2,    &  x_2 x_1&= q_{21} x_1 x_2,   \label{eq:5.2}\\
x_{\frac{3}{2}2}^2    &=0,                           &x_{\frac{3}{2} \frac 52}^{2}&=0,  \label{eq:5.3}\\
x_{1\frac 52}&=x_{\frac{3}{2}2},    &x_{\frac 32 \frac 52}x_2&=q_{12}x_2x_{\frac 32 \frac 52},   \label{eq:5.4}\\
x_{\frac{3}{2}2}x_{\frac 52} &= q_{12}x_{\frac 52} x_{\frac{3}{2}2}+ q_{12} x_2x_{\frac 32 2}, &  x_{\frac 32 2}x_{\frac 32 \frac 52}&= x_{\frac 32 \frac 52} x_{\frac 32 2}.\label{eq:5.5}
\end{align}
A basis of $\toba(\sele_{q}(1,0))$ is given by the monomials
\begin{align}
x_{\frac 52}^{m_1} x_{2}^{m_2}x_{\frac{3}{2}2}^{m_3}x_{\frac{3}{2}\frac 52}^{m_4} x_{\frac{3}{2}}^{m_5} x_1^{m_6},\label{eq:5.7}
\end{align}
where $m_i\in \I_{0,1}$ and $i\in \I_{6}$. Hence $\dim \toba(\sele_{q}(1,0)) = 2^6$.
\end{Theorem}
\begin{proof}
The relations \eqref{eq:5.1}, \eqref{eq:5.2} and $x_{\frac32 2}^2=0$ hold in $\toba(\sele_q(1,0))$, since
\[
V_{g_1}\oplus\Bbbk x_2\simeq\eny_{q_{12}}(1)
\qquad\text{and}\qquad
V_{g_2}\oplus\Bbbk x_1\simeq\eny_{q_{21}}(1)
\]
as braided vector spaces. To verify the remaining relations, we use skew
derivations. The relevant computations are summarized in the following
table:
\[\begin{array}{c|ccc}
 & x_{1\frac52} & x_{\frac32\frac52} & x_{\frac32 2} \\ \hline
\partial_1    & 0 & 0 & 0 \\
\partial_{\frac32}    & 0 & 0 & 0 \\
\partial_2    & x_1 & x_{\frac32}+x_1 & x_1 \\
\partial_{\frac52}    & 0 & x_1 & 0
\end{array}\]
From the previous table and $\partial_{2}(x_{1\frac 52}+x_{\frac 32 2})=x_{1}+x_{1}=0$, we conclude that $x_{1\frac 52}=x_{\frac{3}{2}2}$.
Also, \begin{align*}
    &\partial_{2}(x_{\frac 32 \frac 52}x_2+q_{12}x_{2}x_{\frac 32\frac 52})=x_{\frac 32 2}+x_{1\frac 52}=0,\\
    &\partial_{\frac 52}(x_{\frac 32 \frac 52}x_2+q_{12}x_{2}x_{\frac 32\frac 52})=x_{1}x_2+q_{12}x_2 x_1=0,
\end{align*}
which implies that $x_{\frac 32 \frac 52}x_2=q_{12}x_2x_{\frac 32 \frac 52}$. It follows directly from \eqref{eq:5.2} that $x_{\frac 32 2}x_{1}=q_{21}x_{1}x_{\frac 32 2}$. Hence
\begin{align*}
    &\partial_{2}(x_{\frac 32 \frac 52}^{2})= x_{\frac 32 \frac 52}x_{1}+x_{\frac 32 2}x_{1}+q_{21}x_{1}x_{\frac 32 2}+x_{\frac 32 \frac 52}x_{1}=0,\\
    &\partial_{\frac 52}(x_{\frac 32 \frac 52}^{2})=x_{\frac 32 \frac 52}x_{1}+x_{\frac 32 2}x_{1}+q_{21}x_{1}x_{\frac 32 2}+x_{\frac 32 \frac 52}x_{1}=0,
\end{align*}
and hence $x_{\frac{3}{2} \frac 52}^{2}=0$. The relations in \eqref{eq:5.5} follow from the previous table and using that
\begin{align*}
    &\partial_{2}(x_{\frac{3}{2}2}x_{\frac 52}+ q_{12}x_{\frac 52} x_{\frac{3}{2}2}+ q_{12} x_2x_{\frac 32 2})=x_1 x_{\frac 52}+q_{12}x_{\frac 52}x_1+x_{\frac 32 2}+q_{12}x_2x_1=0,\\
    &\partial_{\frac 52}(x_{\frac{3}{2}2}x_{\frac 52}+ q_{12}x_{\frac 52} x_{\frac{3}{2}2}+ q_{12} x_2x_{\frac 32 2})=x_{\frac 32 2 }+x_{\frac 32 2}=0,\\
    &\partial_{2}(x_{\frac 32 2}x_{\frac 32 \frac 52}+x_{\frac 32 \frac 52}x_{\frac 32 2})=x_{\frac 32 \frac 52}x_{1}+x_{\frac 32 2}x_{1}+q_{21}x_{1}x_{\frac 32 2}+x_{\frac 32 \frac 52}x_{1}=0.
\end{align*}

Let $\toba$ be the algebra defined by the presentation above. By the preceding paragraph, there is a surjective algebra morphism $\pi:\toba\longrightarrow\toba(\sele_q(1,0))$. The defining relations further imply
\begin{equation}\label{eq:5.8}
\begin{aligned}
x_1x_{\frac32 2}&=q_{12}x_{\frac32 2}x_1,&x_1x_{\frac32\frac52}&=q_{12}(x_{\frac32\frac52}+x_{\frac32 2})x_1,\\
x_{\frac32}x_{\frac32 2}&=q_{12}x_{\frac32 2}x_{\frac32},&x_{\frac32}x_{\frac32\frac52}&=q_{12}(x_{\frac32 2}+x_{\frac32\frac52})x_{\frac32},\\
x_{\frac32\frac52}x_{\frac52}&=q_{12}(x_{\frac52}+x_2)x_{\frac32\frac52},&x_2x_{\frac32 2}&=q_{21}x_{\frac32 2}x_2.
\end{aligned}
\end{equation}
Together with the definitions of $x_{\frac32 2}$ and $x_{\frac32\frac52}$, these relations show that the monomials in
\eqref{eq:5.7} span $\toba$, and hence also $\toba(\sele_q(1,0))$.

It remains to prove that the monomials in \eqref{eq:5.7} are linearly independent. Suppose, to the contrary, that there exists a nontrivial
homogeneous linear combination $S$ of these monomials that vanishes, and choose such an $S$ of minimal degree. Since
\begin{align*}
\partial_1\bigl(x_{\frac52}^{m_1}x_2^{m_2}x_{\frac32 2}^{m_3}x_{\frac32\frac52}^{m_4}x_{\frac32}^{m_5}x_1^{m_6}\bigr)&=
\delta_{m_6,1}x_{\frac52}^{m_1}x_2^{m_2}x_{\frac32 2}^{m_3}x_{\frac32\frac52}^{m_4}x_{\frac32}^{m_5},\\
\partial_{\frac32}\bigl(x_{\frac52}^{m_1}x_2^{m_2}x_{\frac32 2}^{m_3}x_{\frac32\frac52}^{m_4}x_{\frac32}^{m_5}\bigr)&=
\delta_{m_5,1}x_{\frac52}^{m_1}x_2^{m_2}x_{\frac32 2}^{m_3}x_{\frac32\frac52}^{m_4},\\
\partial_1\partial_{\frac52}\bigl(x_{\frac52}^{m_1}x_2^{m_2}x_{\frac32 2}^{m_3}x_{\frac32\frac52}^{m_4}\bigr)&=
\delta_{m_4,1}x_{\frac52}^{m_1}x_2^{m_2}x_{\frac32 2}^{m_3},\\
\partial_1\partial_2\bigl(x_{\frac52}^{m_1}x_2^{m_2}x_{\frac32 2}^{m_3}\bigr)&=\delta_{m_3,1}x_{\frac52}^{m_1}x_2^{m_2},
\end{align*}
the minimality of $S$ successively forces $m_6=m_5=m_4=m_3=0$. Hence $S$ reduces to a nontrivial linear combination of the monomials
$x_{\frac52}^{m_1}x_2^{m_2}$, which are linearly independent in $\toba(V_{g_2}\oplus \Bbbk x_1)\simeq \toba(\eny_{q_{21}}(1))$, a contradiction. Therefore, the monomials in \eqref{eq:5.7} are linearly independent, and consequently $\toba\simeq\toba(\sele_q(1,0))$.
\end{proof}

\subsubsection{A realization of $\toba(\sele_q(1,0))$}

Let $q_{12}$ be a primitive $n$-th root of unity, and let
\[\widetilde{\Gamma}=C_m\times C_m, \]
where $m=\operatorname{lcm}(n,2)$. Then $\sele_q(1,0)\in\yd{\Bbbk\widetilde{\Gamma}}$, with action and coaction
as in the case of $\Bbbk\Gamma$ considered in the previous subsection.

\begin{Proposition}
Let $H=\toba(\sele_q(1,0))\#\Bbbk\widetilde{\Gamma}$ and
$\widetilde{\Gamma}=\langle g_1\rangle\times\langle g_2\rangle$.
Then $H$ is a finite-dimensional pointed Hopf algebra generated by
$x_1,x_{\frac32},x_2,x_{\frac52},g_1,g_2$ satisfying
\eqref{eq:5.1}--\eqref{eq:5.5} and the following relations:
\begin{align*}
g_1x_1&=x_1g_1,&g_2x_1&=q_{21}x_1g_2,\\
g_1x_{\frac32}&=x_{\frac32}g_1,&g_2x_{\frac32}&=q_{21}(x_{\frac32}+x_1)g_2,\\
g_1x_2&=q_{12}x_2g_1,&g_2x_2&=x_2g_2,\\
g_1x_{\frac52}&=q_{12}(x_{\frac52}+x_2)g_1,&g_2x_{\frac52}&=x_{\frac52}g_2,\\
g_1^m&=g_2^m=1,&g_1g_2&=g_2g_1.\end{align*}
The comultiplication is determined by $g_i\in G(H)$, $i\in\I_2$,
$x_i\in P_{1,g_i}(H)$, $i\in\I_2$, $x_{\frac32}\in P_{1,g_1}(H)$ and $x_{\frac52}\in P_{1,g_2}(H)$.
A basis of $H$ is given by the monomials
\[x_{\frac52}^{r_1}x_2^{r_2}x_{\frac32 2}^{r_3}x_{\frac32\frac52}^{r_4}x_{\frac32}^{r_5}x_1^{r_6}g_1^{s_1}g_2^{s_2},\]
where $r_i\in\I_{0,1}$, $i\in\I_6$, $s_j\in\I_{0,m-1}$, $j\in\I_2$. Hence $\dim H=2^6m^2$.
\end{Proposition}

\subsection{Nichols algebra of a pale block and a Jordan block} Throughout this subsection, $a=b=1$ and we use the notation
\[\mathtt{t}=x_{\frac{3}{2}\frac 52}x_{\frac 52}+q_{12}x_{\frac 52}x_{\frac{3}{2}\frac 52},\quad \mathtt{w}=x_{\frac{3}{2}\frac 52}x_{2}+q_{12}x_{2}x_{\frac{3}{2}\frac 52}, \quad \mathtt{x}=x_{\frac{3}{2}\frac{5}{2}} x_{\frac{3}{2}2} + x_{\frac{3}{2}2} x_{\frac{3}{2}\frac{5}{2}}.\]
For this case we have the following result.

\begin{Theorem}\label{pale-Jordan}
The Nichols algebra $\toba(\sele_{q}(1,1))$ is generated by 
$x_1$, $x_{\frac{3}{2}}$, $x_2$ and $x_{\frac{5}{2}}$ with defining relations
\begin{align}
x_1^2 &= 0, &\qquad x_{\frac{3}{2}}^2 &= 0, \label{eq:rel1}\\
x_1 x_{\frac{3}{2}} &=  x_{\frac{3}{2}} x_1,  &\qquad x_1 x_2 &= q_{12} x_2 x_1, \label{eq:rel2}\\
x_{\frac{3}{2}2}^2 &= 0, &\qquad x_2 x_{\frac{3}{2}2} &= q_{21} x_{\frac{3}{2}2} x_2,\label{eq:rel3}\\
x_{2}^2&=0, &\qquad x_{\frac 52}^{4}&=0,\label{eq:rel4}\\
x_2x_{\frac 52}x_2x_{\frac 52}&=x_{\frac 52}x_2x_{\frac 52}x_2, &\qquad x_{\frac 52}^2x_{2}&=x_{2}x_{\frac 52}^2+x_2x_{\frac 52}x_2, \label{eq:rel5} \\
x_{1\frac{5}{2}}& = x_{\frac{3}{2}2}, &\qquad x_{\frac 32 \frac 52}\mathtt{w}&=q_{12}\mathtt{w}x_{\frac 32 \frac 52}, \label{eq:rel6}\\
x_{\frac{3}{2}2} x_{\frac{5}{2}} &= q_{12} x_{\frac{5}{2}} x_{\frac{3}{2}2} + \mathtt{w},&\qquad x_{\frac 32 \frac 52}^{4}&=0,\label{eq:rel7}\\
x_{\frac{3}{2}2}\mathtt{t} &= q_{12} (\mathtt{t} +\mathtt{w})\, x_{\frac{3}{2}2},&\qquad \mathtt{w}\, x_{\frac{5}{2}} &= q_{12} (x_{\frac{5}{2}} + x_2)\, \mathtt{w},\label{eq:rel8}\\
x_{\frac{3}{2}\frac 52} \mathtt{t} &= q_{12} (\mathtt{t} +\mathtt{w})\, x_{\frac{3}{2}\frac 52}, &\qquad            \mathtt{t} x_{\frac{5}{2}} &= q_{12} (x_{\frac{5}{2}} + x_2)\, \mathtt{t},    \label{eq:rel9}\\
\mathtt{w}^2&=\mathtt{t}^2=\mathtt{x}^2=0.\label{eq:rel10}
\end{align}
A basis of $\toba(\sele_{q}(1,1))$  is given by the monomials
\begin{align}
x_2^{m_1}x_{2\frac 52}^{m_2}x_{\frac 52}^{n_1}\mathtt{w}^{m_3}\mathtt{t}^{m_4} x_{\frac{3}{2}\frac{5}{2}}^{n_2}  \mathtt{x}^{m_5} x_{\frac{3}{2}2}^{m_6} x_{\frac{3}{2}}^{m_7} x_1^{m_8},
\label{eq:rel11}
\end{align}
where $m_i\in\I_{0,1}$ for $i\in\I_8$ and
$n_j\in\I_{0,3}$ for $j\in\I_2$. Hence $\dim \toba(\sele_{q}(1,1)) = 2^{12}$.
\end{Theorem}

\begin{proof}
We prove that the relations \eqref{eq:rel1}--\eqref{eq:rel10} hold in $\toba(\sele_q(1,1))$. First, the relations \eqref{eq:rel1},
\eqref{eq:rel2} and \eqref{eq:rel3} hold since $V_{g_1}\oplus\Bbbk x_2\simeq\eny_{q_{12}}(1)$ as braided vector spaces. Also \eqref{eq:rel4} and
\eqref{eq:rel5} follow from $\langle x_2,x_{\frac52}\rangle\simeq\mathcal V(1,2)$
as braided vector spaces; see Proposition~\ref{nichols-jordanplane}.

To verify the remaining relations, we use skew derivations. The relevant
computations with skew derivations are summarized in the following table:
\begin{table}[ht]
\centering
\begin{math}
\begin{array}{c|ccccccc}
&x_{\frac32 2}&x_{\frac32\frac52}&x_{2\frac52}&x_{1\frac52}&\mathtt{t}&\mathtt{w}&\mathtt{x}\\ \hline
\partial_1&0&0&0&0&0&0&0\\\partial_{\frac32}&0&0&0&0&0&0&0\\
\partial_2&x_1&x_{\frac32}+x_1&x_2&x_1&x_{\frac32\frac52}&x_{\frac32 2}&x_{\frac32 2}x_1\\
\partial_{\frac52}&0&x_1&0&0&x_{\frac32 2}&0&0\end{array}
\end{math}
\caption{Skew derivations}
\label{tab:skew-pale-Jordan}
\end{table}

By the previous table, $\partial_i(x_{1\frac52}+x_{\frac32 2})=0$, for $i\in \{1,\frac 32,\frac 52\}$. Moreover, $\partial_{2}(x_{1\frac 52}+x_{\frac 32 2})=x_1+x_1=0$ and hence the relation on the left-hand side of \eqref{eq:rel6} is proved. We claim that the following relations hold:
\begin{equation}\label{eq:reln1}
       \begin{aligned}
            x_1x_{\frac32\frac52} &= q_{12}(x_{\frac32\frac52}+x_{\frac32 2})x_1,  &\,\,\, x_{\frac32}x_{\frac32\frac52} &= q_{12}(x_{\frac32\frac52}+x_{\frac32 2})x_{\frac32},\\
            x_1x_{\frac32 2}& = q_{12}x_{\frac32 2}x_1, &\,\,\, \mathtt{x}x_{\frac32 2} &= x_{\frac32 2}\mathtt{x},\\
            x_1\mathtt{w} &= q_{12}^2\mathtt{w}x_1,&\,\,\,  x_{\frac32}\mathtt{w} &= q_{12}\mathtt{x}+q_{12}^2\mathtt{w}x_{\frac32},\\           
            x_1\mathtt{t} &= q_{12}\mathtt{x}+q_{12}^2\mathtt{t}x_1,  & \,\,\, x_{\frac32}x_{2\frac52}& = q_{12}^2x_{2\frac52 }x_{\frac32}+q_{12}x_2x_{\frac32 2},\\
            x_{\frac32}\mathtt{t} &= q_{12}^2\mathtt{t}x_{\frac32}+q_{12}x_{\frac32 2}x_{\frac32\frac52},    &\,\,\,x_{\frac32 2}\mathtt{w}&= q_{12}\mathtt{w}x_{\frac32 2}.
        \end{aligned}
\end{equation}
Indeed, the relations in the first three lines of \eqref{eq:reln1} follow from the relations already established. We shall also use the relation (see \eqref{eq:rel7}) $x_{\frac32 2}x_{\frac52}=q_{12}x_{\frac52}x_{\frac32 2}+\mathtt{w}$, which follows by applying the skew derivations in Table~\ref{tab:skew-pale-Jordan}. The relations in the last two lines of \eqref{eq:reln1} then follow from this relation together with the preceding ones.

Using the relations already established, we obtain $\partial_2\bigl(x_{\frac32\frac52}\mathtt w+q_{12}\mathtt w x_{\frac32\frac52}\bigr)=0$ and $\partial_{\frac52}\bigl(x_{\frac32\frac52}\mathtt w+q_{12}\mathtt w x_{\frac32\frac52}\bigr)=0$. Thus the relation on the right-hand side of \eqref{eq:rel6} follows. In order to prove the relation on the right-hand side of \eqref{eq:rel7}, we use Table~\ref{tab:skew-pale-Jordan} and \eqref{eq:reln1} to conclude that
\begin{flalign*}
    \partial_2(x_{\frac 32 \frac 52}^{2})&=q_{21}(x_1+x_{\frac 32})x_{\frac 32 \frac 52}+ x_{\frac 32 \frac 52}(x_1+x_{\frac 32})=x_{\frac 32 2}(x_1+x_{\frac 32}),\\
   &\partial_{\frac 52}(x_{\frac 32 \frac 52}^{2})=q_{21}x_1x_{\frac 32\frac 52}+x_{\frac 32\frac 52}x_1=x_{\frac 32 2}.
\end{flalign*}
Hence, $\partial_2(x_{\frac 32\frac 52}^{4})=\partial_{\frac 52}(x_{\frac 32\frac 52}^{4})=0$ which implies that $x_{\frac 32\frac 52}^{4}=0$.
Proceeding similarly, we use \eqref{eq:reln1}, Table~\ref{tab:skew-pale-Jordan} and the relations already established to prove \eqref{eq:rel8}, \eqref{eq:rel9} and
\eqref{eq:rel10}.

Let $\toba$ be the algebra defined by the presentation in the statement of the theorem. By the preceding paragraph, there is a surjective algebra morphism $\pi:\toba\longrightarrow\toba(\sele_q(1,1))$. The defining relations further imply
\begin{equation}\label{eq:rel12}
\begin{aligned}
x_1x_{2\frac52}  &= q_{12}^2x_{2\frac52}x_1,&\qquad x_1\mathtt{x}  &= q_{12}^2\mathtt{x}x_1,\\
x_{\frac32}x_{\frac32 2}  &= q_{12}x_{\frac32 2}x_{\frac32},&\qquad x_{\frac32}\mathtt{x}  &= q_{12}^2\mathtt{x}x_{\frac32},\\
x_{\frac32\frac52}x_{2\frac52}  &= q_{12}^2x_{2\frac52}x_{\frac32\frac52} +q_{12}x_2\mathtt{w},&\qquad x_{\frac32\frac52}\mathtt{x}  &= \mathtt{x}     \bigl(x_{\frac32\frac52}+x_{\frac32 2}\bigr),\\
x_{\frac32 2}x_{2\frac52}  &= q_{12}^2x_{2\frac52}x_{\frac32 2},&\qquad \mathtt{x}\mathtt{t}  &= q_{12}^2\mathtt{t}\mathtt{x},\\
\mathtt{x}\mathtt{w}  &= q_{12}^2\mathtt{w}\mathtt{x},&\qquad\mathtt{x}x_{\frac52}  &= q_{12}^2x_{\frac52}\mathtt{x}     +q_{12}\mathtt{w}x_{\frac32 2},\\
\mathtt{x}x_2  &= q_{12}^2x_2\mathtt{x},&\qquad\mathtt{x}x_{2\frac52}  &= q_{12}^4x_{2\frac52}\mathtt{x},\\\mathtt{t}x_2  &= q_{12}x_2\mathtt{t},&\qquad\mathtt{t}x_{2\frac52}  &= q_{12}^2x_{2\frac52}\mathtt{t},\\
\mathtt{t}\mathtt{w}  &= \mathtt{w}\mathtt{t},&\qquad\mathtt{w}x_2  &= q_{12}x_2\mathtt{w},\\
\mathtt{w}x_{2\frac52}  &= q_{12}^2x_{2\frac52}\mathtt{w},&\qquad x_{\frac52}x_{2\frac52}  &= x_{2\frac52}x_{\frac52}+x_2x_{2\frac52},\\
x_{2\frac52}x_2  &= x_2x_{2\frac52}.
\end{aligned}
\end{equation}
Together with the definitions of $\mathtt{t}$, $\mathtt{w}$ and $\mathtt{x}$, relations \eqref{eq:rel1}--\eqref{eq:rel10},
\eqref{eq:reln1} and \eqref{eq:rel12} show that the monomials in \eqref{eq:rel11} span $\toba$, and hence also $\toba(\sele_q(1,1))$.

It remains to prove that the monomials in \eqref{eq:rel11} are linearly independent. Suppose, to the contrary, that there exists a nontrivial homogeneous linear combination $S$ of these monomials that vanishes and choose such an $S$ of minimal degree.
First, for $M=x_2^{m_1}x_{2\frac52}^{m_2}x_{\frac52}^{n_1}\mathtt{w}^{m_3}\mathtt{t}^{m_4}x_{\frac32\frac52}^{n_2}\mathtt{x}^{m_5}
x_{\frac32 2}^{m_6}x_{\frac32}^{m_7}x_1^{m_8}$ we have
\begin{align*}
\partial_1(M)&=\delta_{m_8,1}x_2^{m_1}x_{2\frac52}^{m_2}x_{\frac52}^{n_1}\mathtt{w}^{m_3}\mathtt{t}^{m_4}
x_{\frac32\frac52}^{n_2}\mathtt{x}^{m_5}x_{\frac32 2}^{m_6}x_{\frac32}^{m_7},\\
\partial_{\frac32}(M)&=\delta_{m_7,1}x_2^{m_1}x_{2\frac52}^{m_2}x_{\frac52}^{n_1}\mathtt{w}^{m_3}\mathtt{t}^{m_4}x_{\frac32\frac52}^{n_2}\mathtt{x}^{m_5}
x_{\frac32 2}^{m_6}.
\end{align*}
Hence, the minimality of $S$ forces $m_8=m_7=0$. Set $D:=\partial_1\partial_2$. For the monomial $A=x_2^{m_1}x_{2\frac52}^{m_2}x_{\frac52}^{n_1}\mathtt{w}^{m_3}\mathtt{t}^{m_4}$ the following identities hold:
\begin{align*}D\bigl(A\mathtt{x}^{m_5}x_{\frac32 2}\bigr)&=A\mathtt{x}^{m_5},\\
D\bigl(Ax_{\frac32\frac52}\mathtt{x}^{m_5}x_{\frac32 2}\bigr)&=A\mathtt{x}^{m_5}x_{\frac32 2}+Ax_{\frac32\frac52}\mathtt{x}^{m_5},\\
D\bigl(Ax_{\frac32\frac52}^{2}\mathtt{x}^{m_5}x_{\frac32 2}\bigr)&=Ax_{\frac32\frac52}^{2}\mathtt{x}^{m_5},\\
D\bigl(Ax_{\frac32\frac52}^{3}\mathtt{x}^{m_5}x_{\frac32 2}\bigr)&=(m_5+1)A\mathtt{x}x_{\frac32 2}+Ax_{\frac32\frac52}^{2}\mathtt{x}^{m_5}x_{\frac32 2}
+Ax_{\frac32\frac52}^{3}\mathtt{x}^{m_5}.
\end{align*}
These identities, together with the minimality of $S$, force $m_6=0$. If $m_5=1$, then
\begin{align*}
D^2(A\mathtt{x})&=A, &\,\,\, D^2(Ax_{\frac32\frac52}\mathtt{x})&=Ax_{\frac32\frac52},\\
D^{2}(Ax_{\frac32\frac52}^{2}\mathtt{x})&=A\mathtt{x}+Ax_{\frac32\frac52}^{2},&\,\,\,D^2(Ax_{\frac32\frac52}^{3}\mathtt{x})&=Ax_{\frac32\frac52}\mathtt{x}+A\mathtt{x}x_{\frac32 2}+Ax_{\frac32\frac52}^{3}.
\end{align*}
Thus, the minimality of $S$ also forces $m_5=0$. Next, set $E:=\partial_1\partial_{\frac52}$. For $1\le n_2\le3$, we have
\begin{align*}
E(Ax_{\frac32\frac52})&=A,& E(Ax_{\frac32\frac52}^{2})&=Ax_{\frac32 2},\\
E(Ax_{\frac32\frac52}^{3})&=A\bigl(x_{\frac32 2}x_{\frac32\frac52}+x_{\frac32\frac52}^{2}\bigr).&&
\end{align*}
It follows successively that the coefficients of the terms with
$n_2=3$, $n_2=2$ and $n_2=1$ vanish. Hence $n_2=0$ for every monomial
occurring in $S$. We now eliminate $m_4$ and $m_3$. Since \[\partial_1\partial_2\partial_{\frac52}\bigl(x_2^{m_1}x_{2\frac52}^{m_2}x_{\frac52}^{n_1}\mathtt{w}^{m_3}\mathtt{t}^{m_4}\bigr)
=x_2^{m_1}x_{2\frac52}^{m_2}x_{\frac52}^{n_1}\mathtt{w}^{m_3},\] the minimality of $S$ gives $m_4=0$. Moreover, \[\partial_1(\partial_2)^2\bigl(x_2^{m_1}x_{2\frac52}^{m_2}x_{\frac52}^{n_1}\mathtt{w}^{m_3}\bigr)=\delta_{m_3,1}x_2^{m_1}x_{2\frac52}^{m_2}x_{\frac52}^{n_1},\]
and therefore $m_3=0$. Consequently, $S$ reduces to a nontrivial linear combination of the $x_2^{m_1}x_{2\frac52}^{m_2}x_{\frac52}^{n_1}$. These monomials are linearly independent in $\toba(\langle x_2,x_{\frac52}\rangle)$, a contradiction. Therefore, the
monomials in \eqref{eq:rel11} are linearly independent.
\end{proof}

\subsubsection{A realization of $\toba(\sele_q(1,1))$}

Let $q_{12}$ be a primitive $n$-th root of unity, and let
\[
\widetilde{\Gamma}=C_m\times C_m,
\]
where $m=\operatorname{lcm}(n,2)$. Then $\sele_q(1,1)\in\yd{\Bbbk\widetilde{\Gamma}}$, with action and coaction
as in the case of $\Bbbk\Gamma$ considered above.

\begin{Proposition}
Let $H=\toba(\sele_q(1,1))\#\Bbbk\widetilde{\Gamma}$ and $\widetilde{\Gamma}=\langle g_1\rangle\times\langle g_2\rangle$.
Then $H$ is a finite-dimensional pointed Hopf algebra generated by $x_1,x_{\frac32},x_2,x_{\frac52},g_1,g_2$ satisfying
\eqref{eq:rel1}--\eqref{eq:rel10} and the following relations:
\begin{align*}
g_1x_1&=x_1g_1,&g_2x_1&=q_{21}x_1g_2,\\
g_1x_{\frac32}&=x_{\frac32}g_1,&g_2x_{\frac32}&=q_{21}(x_{\frac32}+x_1)g_2,\\
g_1x_2&=q_{12}x_2g_1,&g_2x_2&=x_2g_2,\\
g_1x_{\frac52}&=q_{12}(x_{\frac52}+x_2)g_1,&g_2x_{\frac52}&=(x_{\frac52}+x_2)g_2,\\
g_1^m&=g_2^m=1,&g_1g_2&=g_2g_1.
\end{align*}
The comultiplication is determined by $g_i\in G(H)$, $i\in\I_2$, $x_i\in P_{1,g_i}(H)$, $i\in\I_2$, $x_{\frac32}\in P_{1,g_1}(H)$ and
$x_{\frac52}\in P_{1,g_2}(H)$.
A basis of $H$ is given by the monomials
\[
x_2^{r_1}x_{2\frac52}^{r_2}x_{\frac52}^{l_1}
\mathtt{w}^{r_3}\mathtt{t}^{r_4}
x_{\frac32\frac52}^{l_2}\mathtt{x}^{r_5}
x_{\frac32 2}^{r_6}x_{\frac32}^{r_7}x_1^{r_8}
g_1^{s_1}g_2^{s_2},
\]
where $r_i\in\I_{0,1}$, $i\in\I_8$, $l_j\in\I_{0,3}$, $j\in\I_2$ and $s_i\in\I_{0,m-1}$, $i\in\I_2$.
Consequently $\dim H=2^{12}m^2$.
\end{Proposition}

\section{A pale block of dimension $2$ and several points}
In this section, we extend to higher dimensions the braided vector space consisting of a pale block and two points considered in Section~\ref{sec-pale-2-points}. Let $\theta>2$ and
\[V=\bigoplus_{k=1}^{\theta}V_{g_k},\]
where $V_{g_1}=\langle x_1,x_{\frac32}\rangle$ is a pale block of eigenvalue $1$  and $V_{g_k}=\langle x_k\rangle$, $k\in\I_{2,\theta}$. We realize $V$ as a Yetter--Drinfeld module over $\Bbbk\Gamma$, where $\Gamma=\mathbb Z^\theta$ is generated by $g_k$, $k\in\I_\theta$. The
action of $\Gamma$ on $V$ is given by
\[
\begin{aligned}
g_1\cdot x_1&=x_1,&g_1\cdot x_{\frac32}&=x_{\frac32},&g_1\cdot x_k&=q_{1k}x_k,\\
g_2\cdot x_1&=q_{21}x_1,&g_2\cdot x_{\frac32}&=q_{21}(x_{\frac32}+x_1),&g_2\cdot x_k&=q_{2k}x_k,\\
g_j\cdot x_1&=q_{j1}x_1,&g_j\cdot x_{\frac32}&=q_{j1}(x_{\frac32}+a_jx_1),&
g_j\cdot x_k&=q_{jk}x_k,
\end{aligned}
\]
for all $k\in\I_{2,\theta}$ and $j\in\I_{3,\theta}$, where $a_j\in\Bbbk$. The coaction is determined by
\[
\delta(x_k)=g_k\otimes x_k,\quad k\in\I_\theta,\qquad\delta(x_{\frac32})=g_1\otimes x_{\frac32}.
\]
We assume that $q_{kk}=1$, $q_{kj}q_{jk}=1$, for all $j,k\in\I_\theta$. If $a_j\in\Bbbk^\times$ for all $j\in\I_{3,\theta}$, we denote the braided vector space defined above by $\eny^\star_{2,\theta-1}$. If $a_j=0$ for all
$j\in\I_{3,\theta}$, we denote the braided vector space by $\eny_{2,\theta-1}$ and the corresponding Nichols algebra by $\toba(\eny_{2,\theta-1})$. We also consider the mixed case. Let $I,J\subseteq\I_{3,\theta}$ be disjoint subsets such that $I\cup J=\I_{3,\theta}$ and assume that
\[
a_i=0,\quad i\in I,\qquad a_j\in\Bbbk^\times,\quad j\in J.
\]
In this case, we denote the corresponding braided vector space by $\eny^\dagger_{2,\theta-1}(I,J)$.
Throughout this section, we use the notation
\[
w_j:=x_{\frac32}x_j+q_{1j}x_jx_{\frac32},\qquad j\in\I_{2,\theta}.
\]

\subsection{Case $\eny^\star_{2,\theta-1}$} Suppose that $a_j\in \Bbbk^\times$, for all $j\in \I_{3,\theta}$. The next result generalizes Theorem \ref{block2_plus_1point}.

\begin{Theorem}\label{block2_plus_theta-point}
The Nichols algebra $\toba(\eny^\star_{2,\theta-1})$ is generated by $x_1,x_{\frac32}$ and $x_k$, $k\in\I_{2,\theta}$, subject to the following relations:
\begin{align}
x_j^2&=0,\quad j\in\I_\theta,&x_{\frac32}^2&=0,\label{eq1-end-20}\\
x_1x_i&=q_{1i}x_ix_1,&x_1x_{\frac32}&=x_{\frac32}x_1,\label{eq2-end-20}\\
x_ix_k&=q_{ik}x_kx_i,\quad i<k,& x_iw_k&=q_{ik}q_{i1}w_kx_i, \label{eq3-end-20}\\
w_k^2&=0,\label{eq4-end-20}
\end{align}
for $i,k\in\I_{2,\theta}$. A basis of $\toba(\eny^\star_{2,\theta-1})$ is given by the monomials
\begin{equation}\label{monomials-end-20}
x_1^{m_1}x_{\frac32}^{m_{\frac32}}x_2^{m_2}w_2^{p_2}x_3^{m_3}w_3^{p_3}\cdots x_\theta^{m_\theta}w_\theta^{p_\theta},
\end{equation}
where $m_1,m_{\frac32},m_k,p_k\in\{0,1\}$ and $k\in\I_{2,\theta}$. Hence $\dim\toba(\eny^\star_{2,\theta-1})=2^{2\theta}$.
\end{Theorem}

\begin{proof}
The relations \eqref{eq1-end-20}, \eqref{eq2-end-20} and \eqref{eq4-end-20} hold in $\toba(\eny^\star_{2,\theta-1})$ since
\[
\langle x_1,x_{\frac32},x_j\rangle
\simeq \eny_{q_{1j}}(1),
\qquad j\in\I_{2,\theta},
\]
as braided vector spaces. We now verify the remaining relations using skew derivations. Let
$j,k\in\I_{2,\theta}$ with $j\ne k$. Since $\partial_{\frac32}(w_k)=0$, we obtain $\partial_{\frac32}\bigl(x_jw_k+q_{jk}q_{j1}w_kx_j\bigr)=0$.
Moreover
\[
\partial_k(w_k)=\begin{cases}
x_1, & k=2,\\
a_kx_1, & k\in\I_{3,\theta},
\end{cases}
\]
and hence $\partial_k\bigl(x_jw_k+q_{jk}q_{j1}w_kx_j\bigr)=0$. Since $g_j\cdot w_k=q_{j1}q_{jk}w_k$, we also have $\partial_j\bigl(x_jw_k+q_{jk}q_{j1}w_kx_j\bigr)=0$. The remaining skew derivatives vanish as well. Thus $x_jw_k=q_{jk}q_{j1}w_kx_j$ for $j\ne k$. Similarly, applying the skew derivations to $x_jx_k+q_{jk}x_kx_j$ gives $x_jx_k=q_{jk}x_kx_j$ for $j<k$. Hence, all the defining relations hold in $\toba(\eny^\star_{2,\theta-1})$.

Let $\toba$ be the algebra defined by the presentation above. By the preceding paragraph, there is a surjective algebra morphism $\pi:\toba\longrightarrow\toba(\eny^\star_{2,\theta-1})$. The defining relations further imply
\begin{equation}\label{eq:conjunto-end-20}
\begin{aligned}
x_1w_k&=q_{1k}w_kx_1,&\qquad x_{\frac32}w_k&=q_{1k}w_kx_{\frac32},\\
x_kw_k&=q_{k1}w_kx_k,&\qquad w_kw_j&=q_{kj}q_{1j}q_{k1}w_jw_k,
\end{aligned}
\end{equation}
for $j,k\in\I_{2,\theta}$ with $j\ne k$ in the last relation.
Together with the definitions of $w_k$, these relations show that the
monomials in \eqref{monomials-end-20} span $\toba$, and hence also
$\toba(\eny^\star_{2,\theta-1})$.

It remains to prove that the monomials in \eqref{monomials-end-20} are linearly independent. We proceed by induction on $\theta$. The case
$\theta=3$ follows from Theorem~\ref{block2_plus_1point}. Assume that the result holds for $\theta-1$. Set $W=V_{g_1}\oplus\cdots\oplus V_{g_{\theta-1}}$. By the splitting technique, there is an isomorphism $\toba(\eny^\star_{2,\theta-1})\simeq\mathcal K\#\toba(W)$,
where $\toba(W)\simeq\toba(\eny^\star_{2,\theta-2})$ and $\mathcal K^1=\ad_c(\toba(W))(V_{g_\theta})$. The defining relations imply that
$\mathcal K^1=\langle x_\theta,w_\theta\rangle$. Moreover,
\[
x_\theta^2=w_\theta^2=0,
\qquad
x_\theta w_\theta=q_{\theta1}w_\theta x_\theta,
\]
so $\mathcal K$ is spanned by $\{1, x_\theta, w_\theta, x_\theta w_\theta\}$. These elements are linearly independent since $\langle x_1,x_{\frac32},x_\theta\rangle \simeq\eny_{q_{1\theta}}(1)$ as braided vector spaces. Therefore, the set 
$\{x_\theta^{m_\theta}w_\theta^{p_\theta} : m_\theta,p_\theta\in\I_{0,1}\}$ is a basis of $\mathcal K$. By the induction hypothesis, the monomials $x_1^{m_1}x_{\frac32}^{m_{\frac32}}x_2^{m_2}w_2^{p_2}\cdots x_{\theta-1}^{m_{\theta-1}}w_{\theta-1}^{p_{\theta-1}}$
form a basis of $\toba(W)$. The vector space decomposition $\toba(\eny^\star_{2,\theta-1})\simeq\mathcal K\otimes\toba(W)$
then shows that the monomials in \eqref{monomials-end-20} are linearly independent. 
\end{proof}

\subsubsection{A realization of $\toba(\eny^\star_{2,\theta-1})$}
Assume that $q_{jk}$ is a primitive $n_{jk}$-th root of unity for $1\le j<k\le\theta$, and let
\[\widetilde{\Gamma}=C_{m_1}\times\cdots\times C_{m_\theta},\]
where $m_1=\operatorname{lcm}(n_{12},\ldots,n_{1\theta})$, $m_j=\operatorname{lcm}\bigl(2,n_{1j},\ldots,n_{j-1,j},n_{j,j+1},\ldots,n_{j\theta}\bigr)$, for $j\in\I_{2,\theta}$. Then $\eny^\star_{2,\theta-1}\in\yd{\Bbbk\widetilde{\Gamma}}$ with action
and coaction as in the case of $\Bbbk\Gamma$ considered above.

\begin{Proposition}
Let $H=\toba(\eny^\star_{2,\theta-1})\#\Bbbk\widetilde{\Gamma}$ and $\widetilde{\Gamma}
=\langle g_1\rangle\times\cdots\times\langle g_\theta\rangle$. Then $H$ is a finite-dimensional pointed Hopf algebra generated by $x_1$,$\,x_{\frac32}$,$\,x_2,\ldots,x_\theta$,\linebreak $g_1,\ldots,g_\theta$ satisfying \eqref{eq1-end-20}--\eqref{eq4-end-20} and the following
relations:
\begin{align*}
g_1x_i&=q_{1i}x_ig_1,&g_1x_{\frac32}&=x_{\frac32}g_1,\\
g_2x_i&=q_{2i}x_ig_2,&g_2x_{\frac32}&=q_{21}(x_{\frac32}+x_1)g_2,\\
g_kx_i&=q_{ki}x_ig_k,&g_kx_{\frac32}&=q_{k1}(x_{\frac32}+a_kx_1)g_k,\\
g_i^{m_i}&=1,& g_ig_j&=g_jg_i,
\end{align*}
for $i,j\in\I_\theta$ and $k\in\I_{3,\theta}$. 
The comultiplication is determined by $g_i\in G(H)$, $i\in\I_\theta$, $x_i\in P_{1,g_i}(H)$, $i\in\I_\theta$ and $x_{\frac32}\in P_{1,g_1}(H)$. A basis of $H$ is given by the monomials
\[
x_1^{r_1}x_{\frac32}^{r_{\frac32}}x_2^{r_2}w_2^{p_2}x_3^{r_3}w_3^{p_3}\cdots x_\theta^{r_\theta}w_\theta^{p_\theta}
g_1^{s_1}\cdots g_\theta^{s_\theta},
\]
where $r_1,r_{\frac32},r_i,p_i\in\I_{0,1}$, $i\in\I_{2,\theta}$ and $s_i\in\I_{0,m_i-1}$, $i\in\I_\theta$. Hence $\dim H=2^{2\theta}m_1\cdots m_\theta$.
\end{Proposition}

\subsection{Case $\eny_{2,\theta-1}$}

Suppose that $a_j=0$ for all $j\in\I_{3,\theta}$. The following result generalizes Theorem~\ref{block-1point-star}.

\begin{Theorem}\label{several-points-case-zero}
The Nichols algebra $\toba(\eny_{2,\theta-1})$ is generated by $x_1,x_{\frac32}$ and $x_k$, $k\in\I_{2,\theta}$, subject to the following
relations:
\begin{align}
x_{\frac 32}x_1 + x_1x_{\frac 32}&=0, &x_{i}^2&=0, &   i&\in\I_2, \label{1eq_2theta}\\
x_{1}x_2 + q_{12}x_2 x_{1}&=0, & w_2^2&=0,  &x_{\frac 32}^2&=0,\label{2eq_2theta}\\
x_{\frac 32} x_j + q_{1j} x_j x_{\frac 32 }&=0, & x_{1}x_j+q_{1j}x_jx_1&=0, & x_{j}^2&=0, \label{3eq_2theta}\\
x_2x_j+ q_{2j}x_jx_2&=0, & x_jx_k+q_{jk}x_kx_j&=0, &j,k\in \I_{3,\theta}, &\quad j<k.\label{4eq_2theta} 
\end{align}
A basis of $\toba(\eny_{2,\theta-1})$ is given by the monomials
\begin{equation}\label{monomials-end5}
x_2^{m_2}w_2^n x_3^{m_3}\cdots x_\theta^{m_\theta}x_{\frac32}^{m_{\theta+1}}x_1^{m_1},
\end{equation}
where $n,m_i\in\I_{0,1}$ and $ i\in\I_{1,\theta+1}$. Hence $\dim\toba(\eny_{2,\theta-1})=2^{\theta+2}$.
\end{Theorem}

\begin{proof}
Since $\langle x_1,x_{\frac32},x_2\rangle\simeq \eny_{q_{12}}(1)$ as braided vector spaces, the relations \eqref{1eq_2theta} and
\eqref{2eq_2theta} hold in $\toba(\eny_{2,\theta-1})$. The remaining relations \eqref{3eq_2theta} and \eqref{4eq_2theta} can be verified using
skew derivations. Thus, all the defining relations hold in
$\toba(\eny_{2,\theta-1})$. 
Let $\toba$ be the algebra defined by the presentation above. By the
preceding paragraph, there is a surjective morphism $\pi:\toba\longrightarrow\toba(\eny_{2,\theta-1})$. The defining relations further imply
\begin{equation}\label{eq:conjunto-end2theta}
\begin{aligned}
x_1w_2&=q_{12}w_2x_1,&x_{\frac32}w_2&=q_{12}w_2x_{\frac32},\\
x_2w_2&=q_{21}w_2x_2,&x_jw_2&=q_{j1}q_{j2}w_2x_j,\quad j\in\I_{3,\theta}.
\end{aligned}
\end{equation}
Together with the definition of $w_2$, these relations show that the monomials in \eqref{monomials-end5} span $\toba$, and hence also
$\toba(\eny_{2,\theta-1})$.

It remains to prove that the monomials in \eqref{monomials-end5} are linearly independent. Suppose, to the contrary, that there exists a
nontrivial homogeneous linear combination $S$ of these monomials that vanishes, and choose such an $S$ of minimal degree.
Notice that
\begin{align*}
\partial_1(x_2^{m_2}w_2^n x_3^{m_3}\cdots x_\theta^{m_\theta}x_{\frac32}^{m_{\theta+1}}x_1^{m_1})&=\delta_{m_1,1}x_2^{m_2}w_2^n x_3^{m_3}\cdots x_\theta^{m_\theta}x_{\frac32}^{m_{\theta+1}},\\
\partial_{\frac32}\bigl(x_2^{m_2}w_2^n x_3^{m_3}\cdots x_\theta^{m_\theta}x_{\frac32}^{m_{\theta+1}}\bigr)&=\delta_{m_{\theta+1},1}
x_2^{m_2}w_2^n x_3^{m_3}\cdots x_\theta^{m_\theta},\\
\partial_k\bigl(x_2^{m_2}w_2^n x_3^{m_3}\cdots x_k^{m_k}\bigr)&=\delta_{m_k,1}x_2^{m_2}w_2^n x_3^{m_3}\cdots x_{k-1}^{m_{k-1}},
\quad k\in\I_{3,\theta},\\
\partial_1\partial_2\bigl(x_2^{m_2}w_2^n\bigr)&=\delta_{n,1}x_2^{m_2}.
\end{align*}
By the minimality of $S$, $m_1=m_{\theta+1}=m_\theta=\cdots=m_3=n=0$. Hence, $S$ reduces to a linear
combination of the monomials $x_2^{m_2}$, $ m_2\in\I_{0,1}$. Since these monomials are linearly independent, we obtain a contradiction.
Therefore, the monomials in \eqref{monomials-end5} are linearly independent. Consequently, $\pi$ is an isomorphism, and the proof is
complete.
\end{proof}

\subsubsection{A realization of $\toba(\eny_{2,\theta-1})$}
Assume that $q_{jk}$ is a primitive $n_{jk}$-th root of unity for $1\le j<k\le\theta$, and let
\[
\widetilde{\Gamma}=C_{m_1}\times\cdots\times C_{m_\theta},
\]
where $m_1=\operatorname{lcm}(n_{12},\ldots,n_{1\theta})$, $m_2=\operatorname{lcm}(2,n_{12},n_{23},\ldots,n_{2\theta})$
and, for $j\in\I_{3,\theta}$, $m_j=\operatorname{lcm}\bigl(n_{1j},\ldots,n_{j-1,j},n_{j,j+1},\ldots,n_{j\theta}\bigr)$. Then $\eny_{2,\theta-1}\in\yd{\Bbbk\widetilde{\Gamma}}$ with action
and coaction as in the case of $\Bbbk\Gamma$ considered above.

\begin{Proposition}
Let $H=\toba(\eny_{2,\theta-1})\#\Bbbk\widetilde{\Gamma}$ and $\widetilde{\Gamma}=\langle g_1\rangle\times\cdots\times\langle g_\theta\rangle$. 
Then $H$ is a finite-dimensional pointed Hopf algebra generated by $x_1$, $x_{\frac32}$, $x_2,\ldots,x_\theta$, $g_1,\ldots,g_\theta$
satisfying \eqref{1eq_2theta}--\eqref{4eq_2theta} and the following
relations:
\begin{align*}
g_1x_i&=q_{1i}x_ig_1,&g_1x_{\frac32}&=x_{\frac32}g_1,\\
g_2x_i&=q_{2i}x_ig_2,&g_2x_{\frac32}&=q_{21}(x_{\frac32}+x_1)g_2,\\
g_kx_i&=q_{ki}x_ig_k,&g_kx_{\frac32}&=q_{k1}x_{\frac32}g_k,\\
g_i^{m_i}&=1, & g_ig_j&=g_jg_i,
\end{align*}
for $i,j\in\I_\theta$ and $k\in\I_{3,\theta}$. The comultiplication is determined by $g_i\in G(H)$, $x_i\in P_{1,g_i}(H)$, $i\in\I_\theta$ and 
$x_{\frac32}\in P_{1,g_1}(H)$. A basis of $H$ is given by the monomials
\[
x_2^{r_2}w_2^n
x_3^{r_3}\cdots x_\theta^{r_\theta}
x_{\frac32}^{r_{\theta+1}}x_1^{r_1}
g_1^{s_1}\cdots g_\theta^{s_\theta},
\]
where $n,r_i\in\I_{0,1}$, $i\in\I_{1,\theta+1}$, $s_i\in\I_{0,m_i-1}$ and $i\in\I_\theta$.
Hence $\dim H=2^{\theta+2}m_1\cdots m_\theta$.
\end{Proposition}

\subsection{Case $\eny^\dagger_{2,\theta-1}(I,J)$}
Let $I,J\subseteq\I_{3,\theta}$ be disjoint, increasingly ordered subsets such that $I\cup J=\I_{3,\theta}$ and assume that
\[
a_i=0,\quad i\in I,
\qquad
a_j\in\Bbbk^\times,\quad j\in J.
\]

\begin{Theorem}\label{several-casezeroandnot}
The Nichols algebra $\toba(\eny^\dagger_{2,\theta-1}(I,J))$ is generated
by $x_1,x_{\frac32}$ and $x_k$, $k\in\I_{2,\theta}$, subject to the
following relations:
\begin{align}
x_k^2&=0,&x_{\frac32}^2&=0,\quad k\in\I_\theta,\label{eqIJ1}\\
x_1x_{\frac32}&=x_{\frac32}x_1,&x_1x_k&=q_{1k}x_kx_1,\quad k\in\I_{2,\theta},\label{eqIJ2}\\
x_2x_i&=q_{2i}x_ix_2, \quad i\in I,&   x_2x_j&=q_{2j}x_jx_2,\quad \ j\in J,\label{eqIJ3}\\
x_jx_k&=q_{jk}x_kx_j,\quad j<k,&w_j^2&=0,\quad j,k\in J,\ \label{eqIJ4}\\
x_jw_k&=q_{jk}q_{j1}w_kx_j,\quad j\ne k,&x_jw_2&=q_{j2}q_{j1}w_2x_j,\quad j,k\in J,\label{eqIJ5}\\
x_ix_l&=q_{il}x_lx_i,\quad i<l,&x_{\frac32}x_i&=q_{1i}x_ix_{\frac32},\quad i,l\in I,\label{eqIJ6}\\
w_2^2&=0,&x_jx_i&=q_{ji}x_ix_j,\quad i\in I,\ j\in J.\label{eqIJ7}
\end{align}
A basis of $\toba(\eny^\dagger_{2,\theta-1}(I,J))$ is given by the
monomials
\begin{equation}\label{monomials-end-IJ}
x_2^{\alpha_2}w_2^{\beta_2}
\prod_{j\in J}
\bigl(x_j^{\alpha_j}w_j^{\beta_j}\bigr)
\prod_{i\in I}x_i^{\alpha_i}
x_{\frac32}^{\alpha_{\frac32}}x_1^{\alpha_1},
\end{equation}
where $\alpha_t, \beta_r\in\I_{0,1}$, $t\in\I_\theta\cup\left\{\frac32\right\}$ and $r\in J\cup\{2\}$.
Hence $\dim\toba(\eny^\dagger_{2,\theta-1}(I,J))=2^{\theta+2+|J|}$.
\end{Theorem}

\begin{proof}
It is enough to observe that $\toba(\eny^{\dagger}_{2,\theta-1}(I,J))\simeq\toba(V_1)\,\underline{\otimes}\,\toba(V_2)$,
where
\[
V_1=\langle x_1,x_{\frac32}\rangle\oplus\bigoplus_{j\in J\cup \{2\} }\langle x_j\rangle,\qquad V_2=\bigoplus_{i\in I}\langle x_i\rangle,
\]
and \(\underline{\otimes}\) denotes the braided tensor product.   
\end{proof}

\subsubsection{A realization of $\toba(\eny^{\dagger}_{2,\theta-1}(I,J))$} 
Assume that $q_{jk}$ is a primitive $n_{jk}$-th root of unity for $1\le j<k\le\theta$ and let
\[
\widetilde{\Gamma}
=C_{m_1}\times\cdots\times C_{m_\theta},
\]
where, for $j\in J\cup\{2\}$, $m_j=\operatorname{lcm}\bigl(2,n_{1j},\ldots,n_{j-1,j},n_{j,j+1},\ldots,n_{j\theta}\bigr)$
and, for $i\in I\cup\{1\}$, $m_i=\operatorname{lcm}\bigl(n_{1i},\ldots,n_{i-1,i},n_{i,i+1},\ldots,n_{i\theta}\bigr)$. Then $\eny^\dagger_{2,\theta-1}(I,J)\in\yd{\Bbbk\widetilde{\Gamma}}$ with action and coaction as in the case of $\Bbbk\Gamma$ considered above.

\begin{Proposition}
Let $H=\toba(\eny^\dagger_{2,\theta-1}(I,J))\#\Bbbk\widetilde{\Gamma}$ and $\widetilde{\Gamma}=\langle g_1\rangle\times\cdots\times\langle g_\theta\rangle$. Then $H$ is a finite-dimensional pointed Hopf algebra generated by $x_1,x_{\frac32},x_2,\ldots,x_\theta,
g_1,\ldots,g_\theta$ satisfying \eqref{eqIJ1}--\eqref{eqIJ7} and the following relations:
\begin{align*}
g_1x_k&=q_{1k}x_kg_1,&g_1x_{\frac32}&=x_{\frac32}g_1,\\
g_2x_k&=q_{2k}x_kg_2,&g_2x_{\frac32}&=q_{21}(x_{\frac32}+x_1)g_2,\\
g_jx_k&=q_{jk}x_kg_j,&g_jx_{\frac32}&=q_{j1}(x_{\frac32}+a_jx_1)g_j,\quad j\in J,\\
g_ix_k&=q_{ik}x_kg_i,&g_ix_{\frac32}&=q_{i1}x_{\frac32}g_i,\quad i\in I,\\
g_i^{m_i}&=1, & g_ig_j&=g_jg_i,\quad i,j\in\I_\theta,
\end{align*}
for $k\in\I_\theta$. The comultiplication is determined by $g_i\in G(H)$, $x_i\in P_{1,g_i}(H)$, $i\in\I_\theta$, and $x_{\frac32}\in P_{1,g_1}(H)$. A basis of $H$ is given by the monomials
\[
x_2^{\alpha_2}w_2^{\beta_2}
\prod_{j\in J}
\bigl(x_j^{\alpha_j}w_j^{\beta_j}\bigr)
\prod_{i\in I}x_i^{\alpha_i}
x_{\frac32}^{\alpha_{\frac32}}x_1^{\alpha_1}
g_1^{s_1}\cdots g_\theta^{s_\theta},
\]
where $\alpha_t, \beta_r\in\I_{0,1}$, $t\in\I_\theta\cup\left\{\frac32\right\}$, $r\in J\cup\{2\}$ and $s_i\in\I_{0,m_i-1}$. Hence $\dim H=2^{\theta+2+|J|}m_1\cdots m_\theta$.
\end{Proposition}

\medbreak
\noindent \emph{Acknowledgements.} We would like to thank I. Angiono for fruitful discussions.

\end{document}